\documentclass[smallextended,nospthms]{svjour3}

\usepackage{amsfonts,amsmath, amssymb, latexsym,amsrefs,appendix, amsthm}
\usepackage{amscd}
\usepackage[all,cmtip]{xy}
\usepackage{float}
\usepackage{enumerate}
\usepackage{pgf}
\usepackage{tikz}
\usetikzlibrary{arrows}
\newcommand{\midarrow}{\tikz \draw[-triangle 90] (0,0) -- +(.03,0);}

\newcommand{\cat}{\textbf}
\newcommand{\co}{\colon\thinspace}
\newcommand{\coker}{\operatorname{coker}}
\newcommand{\cS}{{\bf S}}

\newcommand{\NN}{{\mathbb N}}

\newcommand{\ZZ}{{\mathbb Z}}
\newcommand{\QQ}{{\mathbb Q}}
\newcommand{\FF}{{\mathbb F}}
\newcommand{\Fp}{{{\mathbb F}_{p}}}

\newcommand{\nod}{\noindent}
\newcommand{\ff}{\infty}
\newcommand{\pt}{\{ \ast \}}
\newcommand{\sset}{{\subseteq}}

\newcommand{\pcp}{{^{\wedge}_p}}

\newcommand{\ratcp}{{^{\wedge}_{\mathbb Q}}}

\newcommand{\hcl}{\text{\textnormal{hocolim}}}

\newcommand{\hl}{\text{\textnormal{holim}}}
\newcommand{\hcIS}{\text{\textnormal{hocolim}}_{I\in {\bf S}}}

\newcommand{\hlIS}{\text{\textnormal{holim}}_{I\in {\bf S}^{op}}}
\newcommand{\map}{\text{\textnormal{map}}}

\newcommand{\im}{\text{\textnormal{Im}}}
\newcommand{\cd}{\text{\textnormal{cd}}}

\newcommand{\BKp}{{BK^{\wedge}_p}}

\newcommand{\rank}{\mathrm{rank}}
\newcommand{\size}{\mathrm{size}}

\newcommand{\Y}{\Gamma}

\newcommand{\Nyk}{N^\Y_K}
\newcommand{\Nyprimek}{N^{\Y'}_K}
\newcommand{\Vyk}{V^\Y_K}

\newcommand{\Vtk}{V^T_K}

\newcommand{\beq}{\begin{eqnarray}}
\newcommand{\eeq}{\end{eqnarray}}
\newcommand{\ari}{\ar@{^{(}->}}
\newcommand{\arif}{\ar@{_{(}->}}
\newcommand{\arline}{\ar@{-}}
\newcommand{\arele}{\ar@{|->}}

\newtheorem{introtheorem}{Theorem}

\newtheorem{thm}{Theorem}[section]
\newtheorem{lem}[thm]{Lemma}
\newtheorem{prop}[thm]{Proposition}

\theoremstyle{definition}
\newtheorem*{ack}{Acknowledgements}
\newtheorem{defn}[thm]{Definition}

\newtheorem{exmp}[thm]{Example}
\newtheorem{remk}[thm]{Remark}
\newtheorem{quest}[thm]{Question}

\theoremstyle{remark}
\newtheorem*{nota}{Conventions, etc}

\begin{document}

\title{Recognizing nullhomotopic maps into the classifying space of a Kac--Moody group\thanks{Supported by the Danish National Research Foundation through the Centre for Symmetry and Deformation (DNRF92).}}

\author{John D. Foley}

\institute{
Department of Mathematical Sciences, University of Copenhagen, Universitetsparken 5,\\ DK-2100 K{\o}benhavn {\O} \\
\email{foleyjd@gmail.com}
}

\titlerunning{Nullhomotopic maps into a Kac--Moody classifying space}
%[Nullhomotopic maps into a Kac--Moody classifying space]
%[J. D. Foley]

\date{Received: date / Accepted: date}
%\date{\today}

\maketitle

\begin{abstract}

This paper extends certain characterizations of nullhomotopic maps between $p$-compact groups
to maps with target the $p$-completed classifying space of a connected Kac--Moody group
and source the classifying space of either a $p$-compact group or  a connected Kac--Moody group.
A well known inductive principle % (and inductive arguments)
 for $p$-compact groups is applied to obtain general, mapping space level results.
An arithmetic fiber square computation shows that a null map from the classifying space of a connected compact Lie group to the classifying space of a connected topological Kac--Moody group
can be detected  by restricting to the maximal torus.
 Null maps between the classifying spaces of connected topological Kac--Moody groups
cannot, in general,   be detected  by restricting to the maximal torus due
to the nonvanishing of an explicit abelian group of obstructions described here.
Nevertheless, partial results are obtained via the application of algebraic discrete Morse theory to higher derived limit calculations. These partial results  %and invariant theory
show that null maps are detected  by restricting to the maximal torus in many cases of interest.
%06A07, 57M07, 22E65
\end{abstract}
%\maketitle

\section{Introduction}

After the proof of the 
Sullivan conjecture by Carlson, Lannes and Miller \cite{Millersulconj, JLsulconj, GCsulconj}, determining the homotopy type of the mapping space $\map(BG', BG)$ for $G'$ and $G$ compact Lie groups became much more technically accessible.
Among the first and most broadly applicable results where characterizations of
nullhomotopic maps by  Friedlander--Mislin \cite{FM2}, Zabrodsky \cite{Zmaps}, and Jackowski--McClure--Oliver {\cite{JMOself}}.
Later, M{\o}ller \cite{JMratiso} applied an inductive principle \cite[9.2]{centers} % (and inductive arguments)
to obtain characterizations of nullhomotopic maps between $p$-compact groups.  Kac--Moody groups provide a generalization of compact Lie groups and
characterizations of null maps between the classifying spaces of rank two (and non-affine) Kac--Moody groups have been given by Aguad\'e--Ruiz \cite{BKtoBK}.

This article adapts the inductive arguments of \cites{centers, JMratiso} to recognize null maps between Kac--Moody group classifying spaces of all ranks.
Throughout the term Kac--Moody group refers to the
unitary form associated to some Kac--Moody Lie algebra \cite[7.4]{kumar} (consistent with \cites{nitutkm, rank2mv, BKtoBK, classKM}).  %By the time this writing, Kac--Moody groups are well-established objects of studyz
We refer the reader to \cite{kumar} for a thorough account of Kac--Moody groups and  \cites{nitutkm} for an introduction to their topology.  See Appendix \ref{append}
for a discussion of how our results apply to Kac--Moody groups associated to derived Kac--Moody Lie algebras. % \cite[7.4]{kumar}

All the previous work characterizing  mapping spaces mentioned above relies on $p$-local techniques. Here we obtain
the following
 characterization of nullhomotopic maps with target the $p$-completed classifying space of a Kac--Moody group.
  \begin{introtheorem}%[cf. {\cite[3.11]{JMOself}}, {\cite[8.1]{BKtoBK}}]
If $K$ is a Kac--Moody group and $\Y$ is a connected $p$-compact or
 Kac--Moody group with maximal torus ${\iota}\co T {\rightarrow} \Y$,
 then the null component of the unpointed mapping space $\map(B\Y, \BKp)_0$ is canonically equivalent to $\BKp$
 and ${f}\co B\Y {\rightarrow} \BKp$ is nullhomotopic  if and only if
$f|_{BT}$ is nullhomotopic.
\label{intronullconnectedBKp}
\end{introtheorem}
\nod
Under the more restrictive hypothesis that $\Y$ and $K$ are rank two (and non-affine) Kac--Moody groups, Theorem \ref{intronullconnectedBKp} is shown in \cite{BKtoBK} but not stated explicitly;
see \cite[Theorem 3.11]{JMOself} for a compact Lie group statement.

The problem of recognizing when integral maps $B\Y \stackrel{g}{\rightarrow} BK$, where now $\Y$ is a connected compact Lie or
 Kac--Moody group, are nullhomotopic exposes new subtleties.
For the cases studied in \cite{BKtoBK, JMOself}, it is proved therein that the homotopy class of $g$ is detected by the collection of homotopy classes of its $p$-completions $g\pcp$ for all primes $p$.
Hence, Theorem \ref{intronullconnectedBKp} implies that null $g$ are detected by restricting to the maximal torus.
% Hence, recognizing nullhomotopic $g$ reduces to recognizing null $g\pcp$.
In our more general setting, we prove that there exist homotopically nontrivial $B\Y \stackrel{g}{\rightarrow} BK$ such that all $g\pcp$ {\em and} the rationalization $g^{{\wedge}}_\QQ$ are null, so
Theorem \ref{intronullconnectedBKp} alone cannot recognize  nullhomotopic $g$.

Nevertheless, we show that
 $g|_{BT}$ is null still implies that $g^{{\wedge}}_\QQ$ is null (see Sections \ref{sec:ratnullkacmoody}-\ref{sec:vanish}) and:

  \begin{introtheorem}
If $K$ is a Kac--Moody group and $\Y=G$ is a connected compact Lie group,
 then ${g}\co BG {\rightarrow} BK$ is nullhomotopic  if and only if
$g|_{BT}$ is nullhomotopic.
\label{nullconnectedLieintro}
\end{introtheorem}
\nod
However---if $\Y$ is not compact---then
 explicit obstructions parameterize the homotopy classes of  $B\Y \stackrel{g}{\rightarrow} BK$ such that $g|_{BT}$ is null
(see Theorem  \ref{nullconnectedcorr}). These obstructions are frequently nonzero,  e.g. if $\Y=K$ is affine or $K$ is compact Lie but $\Y$ is not
(see Examples \ref{exmp:KA} and \ref{1spherical}, respectively).

 These negative results  may disappoint some readers, % focused on affine Kac--Moody groups,
  but they do not preclude a Kac--Moody analog of the %now classic
description given in \cite[Theorems 2 and 3]{JMOself} of
$\map(BG, BG)$  for $G$ a compact {\em simple} Lie group. Because the standard representation of the Weyl group
of an affine Kac--Moody group is not irreducible, affine $K$ are not natural analogs of simple $G$.
In contrast, we prove a partial step toward such a result.
  \begin{introtheorem}
If $K$ is an indecomposable, non-affine and 2-spherical Kac--Moody group,
  then
${g}\co BK {\rightarrow} BK$ is nullhomotopic if and only if $g|_{BT}$ is nullhomotopic.
\label{nullconnectedintro}
\end{introtheorem}
\nod
The restriction to indecomposable $K$ avoids direct product decompositions 
while the 2-spherical condition (see Definition \ref{def:nspherical}) implies a nilpotency  hypothesis on $H^\ast(BK,\QQ)$ that is necessary for the conclusion that ${g}$ is nullhomotopic if and only if $g|_{BT}$ is nullhomotopic.
In particular,  we verify this nilpotency hypothesis by direct computation and give examples of $H^\ast(BK,\QQ)$ as a ring  which have not appeared previously (see Examples \ref{1spherical} and \ref{2sphericalcont}). Moreover,
 we do not know of any of example of an indecomposable, non-affine Kac--Moody group for which this nilpotency  hypothesis fails
(see Question \ref{quest:limvanish}). These computations are facilitated---technically and conceptually---by the application of algebraic discrete Morse theory \cites{aDMT,aDMTlong,CGNDMT}.
Theorem  \ref{nullconnectedintro} can be interpreted as a uniqueness result
for the trivial unstable Adams operation (cf. \cite[Theorem 1]{JMOself} and \cite[$\S$ 10]{BKtoBK}).  Though unstable Adams operations for Kac--Moody groups have been constructed,
their homotopical uniqueness %of these maps
is not settled for rank greater than 2
  \cite[Theorem D and Question 4.4]{PhDart}).

 Given the success found by studying $p$-compact groups and extending Lie theory to Kac--Moody theory (see, e.g., \cites{homotopyLie, JG} and \cites{kumar}, resp.),
  it is hoped that aspects of $p$-compact group theory will extend to a ``homotopy version of a Kac--Moody group'' at $p$
called a homotopy Kac--Moody group.
This paper is uncommonly  successful in repurposing
$p$-compact group methodology (cf. \cites{nitutkm, nituthesis, rank2mv, BKtoBK}), especially inductive techniques;
see Theorem \ref{nullBX'connected} for a generalization of Theorem \ref{intronullconnectedBKp} to this perspective.
% The present paper supports this perspective by demonstrating that some fundamental structure results
%from \cite{JMratiso}
%can be extended to a broad class of spaces that are homotopically similar to the $p$-completed classifying space of a Kac--Moody group.
%
For instance,
 the class of spaces to which our $p$-local results apply includes all representatives  of the homotopy types $\mathcal{S}^\ast$
 introduced in \cites{rank2mv}.  Elements of  $\mathcal{S}^\ast$ are constructed as pushouts of $p$-compact groups and generalize $p$-completed, rank two, Kac--Moody group classifying spaces.
 %
% Indeed, our structure theorems are put into practice
% in a forthcoming paper \cite{maps} where maps between representatives  of $\mathcal{S}^\ast$ are studied.

While the analogy to compact Lie groups is strong at $p$,
this paper gives explicit evidence that the classifying spaces of Kac--Moody groups are rationally less analogous  to the classifying spaces of compact Lie groups.
It has already been observed that constructing integral mapping spaces via the arithmetic fibre square results in new issues---even if $BK$ is rationally an Eilenberg--MacLane space (cf. \cite[Remark 3.4, Theorem 6.2]{BKtoBK}).
We believe that the problems of describing the rational homotopy type of  and constructing integral mapping spaces between
the classifying spaces of Kac--Moody groups
deserve  investigation apart from $p$-local concerns.
We hope that
our more preliminary rational and integral results stimulate  further research in these directions. Our $H^\ast(BK,\QQ)$  computations may find independent interest
 and efforts have been made so that they can be read somewhat independently. %, according to the reader's preference.

 \begin{nota}
Since the author wishes to highlight the reuse of $p$-compact group methods,
this paper does not attempt to be completely self-contained and points of comparison
are documented throughout.  We follow the (mostly) standard definitions and notations for $p$-compact groups
appearing in \cite{homotopyLie}.  Notably, we refer to  $BY$ as a $p$-compact group
rather than the classifying space of a $p$-compact group $Y:=\Omega BY$.
For the reader's convenience,
we %sometimes provide reminders of definitions,
maintain the following notational conventions:

$K=$ a (minimal, split, unitary) Kac--Moody group \cites{nitutkm}, \cite[7.4]{kumar}

$G=$ a compact Lie group

$\Y=$ a connected compact Lie or Kac--Moody group

$\rank(\Y)=$ the rank of a maximal torus of $\Y$

$BY=$ a $p$-compact group \cite[1.10]{homotopyLie}

$BP=$ a $p$-compact toral group \cite[2.1]{homotopyLie} or, more briefly, a $p$-toral group

$\cd_{\Fp}{Y}=$ the maximal nonzero dimension of $H^\ast( Y, \Fp)$

$(-)_\Fp=$ the $\Fp$-homology localization functor \cite{Fploc}

$(-)_\pcp=$ the Bousfield-Kan $p$-completion (homotopy) functor \cite{BK}

$(-)_\ratcp=$ the Bousfield-Kan rationalization (homotopy) functor \cite{BK}

$BX=$ an analog of $\BKp$ with associated $p$-compact groups $BX_I$

 \nod and recommend having  \cites{centers, JMratiso} at hand while reading Section \ref{sec:reduce}.
Whenever a space level $p$-completion is functor needed, we will be able to apply $(-)_\Fp$ which, for
the spaces under consideration, will coincide with $(-)_\pcp$, up to homotopy; see, e.g., (\ref{hocoBKlocal}) below.

The author can also recommend subsets of the article which may appeal to various readers.
 Sections \ref{sec:induct}-\ref{sec:nullkacmoody}
prove Theorem \ref{intronullconnectedBKp} and constitute the $p$-local, homotopy Kac--Moody group story.  Sections \ref{sec:ratnullkacmoody}-\ref{sec:invarianttheory}  can be read by taking Theorem \ref{intronullconnectedBKp}
as a starting point and place Theorem  \ref{nullconnectedintro} in context, with proof.
Section \ref{sec:ratnullkacmoody} outlines the local to global construction of integral $B\Y \stackrel{g}{\rightarrow} BK$ and includes a proof of Theorem \ref{nullconnectedLieintro}.
It can be read alone or supplemented by Section \ref{sec:vanish} which gives further detail on vanishing and homotopical uniqueness.
Section \ref{sec:higherlimits} displays examples of $H^\ast(BK,\QQ)$
 and introduces an approach to compute $H^\ast(BK,\QQ)$ via algebraic discrete Morse theory.
A reading of Section \ref{sec:higherlimits} can also be supplemented by Section \ref{sec:invarianttheory} where the key vanishing result (Theorem \ref{lim2vanish}) for Theorem  \ref{nullconnectedintro} is proven by inputting invariant theory into
this approach.
 \end{nota}

%Key references: \cites{centers, JMOself, BLOJAMS, JMratiso}
%

\section{Structure results by induction}
\label{sec:induct}
\label{sec:reduce}
\label{sec:ptoral}

This section extends some structural, mapping space level
 results ({\cite[10.1]{centers}} and {\cite[2.11,6.1,6.6,6.7]{JMratiso}}) concerning maps between $p$-compact
groups to maps $BY {\rightarrow} BX$ where: 1) $BY$ is a $p$-compact group; and 2) $BX$ is within a class of potential (connected) homotopy Kac--Moody groups.
Specifically,
we construct the spaces $BX$ in this class as follows.
  \beq
BX :=(\hcIS  B X_I)\pcp = (\hcl_\cS  D)\pcp
\label{BXdefined}
\eeq
\nod where $D\co\cS \rightarrow p\mbox{-}\cat{cpt}$ is some diagram of monomorphisms of connected $p$-compact groups
(i.e. $B X_I=D(I)$ is simply-connected)
such that the index category $\cS$ is a finite, contractible poset.

 % these results,
 As a motivating example, let us recall
how a $p$-completed Kac--Moody group classifying space can be constructed as in (\ref{BXdefined}).
See Section \ref{sec:nullkacmoody} for further examples. Any Kac--Moody group $K$ has
a homotopy decomposition \cite{nitutkm}
\beq
\hcl_{I\in \cat{S}} BK_I \stackrel{\thicksim}{\longrightarrow} BK
\label{hocoBK}
%\nonumber
\eeq
\nod where: 1) $K_I$ are connected Lie groups; (2) $K_I \le K$;  and (3) $\cat{S}$ is the finite poset of spherical subsets  $I \subset \{ 1, \ldots , \size(A)\}$ ordered by inclusion  for $A$ the associated generalized Cartan matrix.

Let $D_K$ denote this  diagram $D_K\co\cS \rightarrow \cat{Spaces}$ given by $I \mapsto BK_I$.
In particular,
$\emptyset \in \cat{S}$ is an initial object so that $|\cat{S}|$ is contractible and
 $BK_\emptyset=BT$ is the classifying space of the maximal torus of $K$.
Moreover, each $I \le J$ in $\cat{S}$ corresponds to an inclusion of subgroups $B(K_I \le K_J)$ and each $\emptyset \le J$ corresponds
 to the
 standard maximal torus inclusion $BT \stackrel{\iota_I}{\rightarrow}BK_I$.
 Since all the spaces involved are simply-connected,
 $p$-completion coincides, up to homotopy, with $\Fp$-homology localization.  Thus,
 we have weak equivalences
\beq
[\hcl_{I\in \cat{S}} (BK_I)_\Fp]_\Fp \stackrel{\thicksim}{\longleftarrow}  [\hcl_{I\in \cat{S}} BK_I]_\Fp \stackrel{\thicksim}{\longrightarrow} BK_\Fp \simeq BK\pcp.
\label{hocoBKlocal}
%\nonumber
\eeq
In particular, $\hcIS  (BK_I)_\Fp$ is simply-connected by Seifert--van Kampen theory (see, e.g., \cite{fun}).
Any $BX$ as in (\ref{BXdefined}) is simply-connected, $p$-complete and $\Fp$-local for the same reasons.
% Since $\emptyset \in \cat{S}$ is an initial object, $|\cat{S}|$ is contractible.

For this section, the key point is to apply an inductive principle to the source $p$-compact group $BY$:
  \begin{thm}[{\cite[9.2]{centers}}]
The only class $\cat{C}$ of $p$-compact groups $BY$ with the five properties:
\begin{enumerate}[(i)]
\item $\cat{C}$ is closed under homotopy equivalence,
\item $\cat{C}$ contains the point,
\item $BY_0 \in \cat{C}$ implies $BY \in \cat{C}$,
\item for simply-connected $BY$, $B(Y/Z) \in \cat{C}$ implies $BY \in \cat{C}$, and
\item for simply-connected $BY \simeq B(Y/Z)$, if $\cd_{\Fp}{Y'}<\cd_{\Fp}{Y}$ implies $BY' \in \cat{C}$, then $BY \in \cat{C}$.
\end{enumerate}
\nod is the class of all $p$-compact groups.
\label{saturated}
\end{thm}
\nod %to the source $p$-compact group $BY$.
In the above, $BY_0$ is the homotopy fiber of the canonical $BY \rightarrow B\pi_1(BY)$
and $ BZ \rightarrow BY \rightarrow B(Y/Z)$ is also a homotopy fiber sequence of $p$-compact groups
where $BZ$ is the identity component of $\map(BY,BY)$.
Inductive arguments will be lifted from \cite[\S 9,10]{centers} and \cite[\S 2,6]{JMratiso}
%as outlined below
to reduce our structure results for $BY {\rightarrow} BX$ to the situation where $BY$ is $p$-toral.
%Proofs for a $p$-toral source are given in Section \ref{sec:ptoral}. 
 See Remark \ref{notremark} for a discussion of alternative approaches.

Let us now state our generalization of \cite[10.1]{centers}.

\begin{thm}
For $BX$ as in (\ref{BXdefined}) the natural map
  \beq
\map(\ast, BX) \longrightarrow \map(BY, BX)_0
\label{nullcentralizersBX}
\eeq
is a weak equivalence for any $p$-compact group $BY$.
\label{nullBX}
\end{thm}

As in \cite{centers}, the following lemma can be applied
to show that it is sufficient to prove
Theorem \ref{nullBX} for $BY$ the classifying space of $p$-toral group.  %whose proof is postponed to \ref{sec:ptoral}.

\begin{lem}
 Let $X$ be a $\Fp$-local space %and $BY$ be a $p$-compact group.
such that the canonical map
  \beq
\map(\ast, X) \longrightarrow \map(BY, X)_0
\label{nullcentralizers}
\eeq
is a weak homotopy equivalence for all $BY$  a $p$-toral $p$-compact group.
Then (\ref{nullcentralizers}) is a weak equivalence for any $p$-compact group $BY$.
\label{null}
\end{lem}

\nod{\bf Proof:} Following \cite[\S 9,10]{centers}, it enough to show that the class of $p$-compact groups for which (\ref{nullcentralizers}) is a weak equivalence
 has  the five properties of Theorem \ref{saturated}. Clearly, (i) and (ii) hold.

  Let us now check (iii).
 By our inductive assumption,
the canonical
   \beq
\map(\ast, X) \longrightarrow \map(BY_0, X)_0
\label{nullcentralizersnot}
\eeq
\nod is  a weak homotopy equivalence.   Setting $\pi:=\pi_1(BY)$ and noting $BY \simeq (BY_0)_{h\pi}$, we obtain a weak equivalence
  \beq
\map(B\pi, X) \longrightarrow \coprod_{\phi|{BY_0}\simeq 0} \map(BY, X)_\phi
\nonumber
\label{nullcentralizerspi}
\eeq
which restricts to a weak equivalence of components
  \beq
\map(B\pi, X)_0 \longrightarrow \map(BY, X)_0 .
\label{nullcentralizerspi}
\eeq
By hypothesis, (\ref{nullcentralizers}) with $BY$ replaced by $B\pi$---a $p$-toral group---is also a weak equivalence. Composing with (\ref{nullcentralizerspi})  yields that  $BY \rightarrow B\pi \rightarrow \ast$
induces the desired weak equivalence.

A similar argument reduces (iv) to the hypothesis that (\ref{nullcentralizers}) is a weak equivalence for $BY=BZ$.  The observation that $BZ$ is a
$p$-toral $p$-compact group \cite{centers} completes the check of (iv).

To check (v), use the homology decomposition for $BY$
  \beq
\hcl_{(V,BV \stackrel{v}{\rightarrow}BY )\in \cat{A}^{op}} \map(BV, BY)_v \longrightarrow BY.
%\nonumber
\label{nullcentralizersdecompostion}
\eeq
\nod where objects of $\cat{A}^{op}$ consist of pairs $(V,BV \stackrel{v}{\rightarrow} BY )$ for $V$ a finite
elementary abelian $p$-group \cite{centers}.
 Now we have weak equivalences
\beq
\map( \ast, X) \rightarrow \map( |\cat{A}^{op}|, X) &\rightarrow& \hl_{\cat{A}}\map( \ast, X) \nonumber \\
 &\rightarrow& \hl_{\cat{A}} \map(\map(BV, BY)_f, X)_0
\nonumber
\label{Alocal}
\eeq
\nod using the fact that the $\Fp$-localization of $|\cat{A}^{op}|$ is contractible  (cf. \cites{centers,JMratiso}).
%In particular, it here that we use the hypothesis that $X$ is $\Fp$-local.
\qed

 By design, the hypotheses of Lemma \ref{null} are compatible with Theorem \ref{nullBX}.

 \nod{\bf Reduction of the proof of Theorem \ref{nullBX}:} %First note that $\hcIS  B X_I$ is simply-connected by Seifert--van Kampen theory (see, e.g., \cite{fun}).
Recalling that $BX$ is $\Fp$-local, we  only need to show that (\ref{nullcentralizersBX})
 is a weak equivalence for
$p$-toral  $BY:=BP$  by Lemma \ref{null}.
\qed

Let us now begin the inductive process  by considering  maps from $B\pi$ into $BX$
for $\pi$ a finite $p$-group.
The following is a specialization of a result of Broto, Levi, and Oliver  \cite[Proposition 4.2]{BLOJAMS}.
Compare \cite[Theorem 6.11]{classKM}.

\begin{prop}
 Let $D\co\cS \rightarrow p\mbox{-}\cat{cpt}$  be a diagram of $p$-compact groups  indexed over a finite poset $\cS$.
For any finite $p$-group $\pi$ the canonical map
 \beq
 [\hcIS \map(B\pi , BX_I)]\pcp  \stackrel{\sim}{\rightarrow} \map(B\pi, (\hcIS  B X_I)\pcp)
  \label{BLOJAMSeq}
\eeq
\nod is a weak homotopy equivalence.
\label{pgrouphocolim}
\end{prop}
\nod{\bf Proof:} To apply \cite[Proposition 4.2]{BLOJAMS}, we must check that $\map(B\pi , BX_I)$ is $p$-complete
and that certain higher limits over $\cS$ vanish in degrees greater than some uniform $n$ depending on $\cS$.
See, for example, \cite[Theorem 5.1]{homotopyLie} for the former condition. For the latter, recall \cite[Appendix II.3]{GZ}, \cite[8.3]{Weibel} that all higher limits over $\cS$
in degrees greater than the maximum chain length of $\cS$ vanish.
\qed

We now apply Proposition \ref{pgrouphocolim} to complete the proof of Theorem \ref{nullBX}.

 \nod{\bf Proof of Theorem \ref{nullBX}:} It remains to verify that (\ref{nullcentralizersBX})
 is a weak equivalence for $BP:=BY$ a
$p$-toral group.
The
existence of a homotopy colimit approximation $\hcl_{n\rightarrow \ff} BP_n \rightarrow BP$ \cite[2.1]{homotopyLie}
 reduces the problem to showing
  \beq
\map(\ast, BX) \longrightarrow \map(BP_n, BX)_0
\label{nullcentralizersBP}
\eeq
is a weak equivalence. By Proposition \ref{pgrouphocolim},
 $\map(BP_n, BX)_0$ is some component of $\hcIS \map(BP_n, BX_I)$.
 In particular,
   \beq
\map(\ast, BX_I) \longrightarrow \map(BP_n, BX_I)_0
\label{nullcentralizersBPI}
\nonumber
\eeq
is a weak equivalence for  each $I \in \cS$ by \cite[9.3]{centers}.
For all $J \le I$ in $\cS$, the space $X_I/X_J$ in the homotopy fiber sequence
   \beq
X_I/X_J \longrightarrow BX_J \longrightarrow BX_I
\label{sullivan}
\eeq
is $\Fp$-finite since $BX_J \rightarrow BX_I$ is a monomorphism. In the induced  homotopy fiber sequence,
   \beq
\map(BP_n ,X_I/X_J) \longrightarrow \map(BP_n,BX_J) \longrightarrow \map(BP_n,BX_I)
\label{sullivan2}
\nonumber
\eeq
\nod the Sullivan conjecture \cite{Millersulconj} gives that
 $\map(BP_n, X_I/X_J) \simeq X_I/X_J$.
Since by construction %the null component
$\map(BP_n, BX_I)_0\simeq BX_I$ is simply-connected, the associated homotopy fiber sequence over $\map(BP_n, BX_I)_0$
implies that
 the set of components of $\map(BP_n, BX_J)$ that push forward into $\map(BP_n, BX_I)_0$ is isomorphic
  $\pi_0(X_I/X_J)$ and  (\ref{sullivan}) implies $\pi_0(X_I/X_J)\cong \pi_0(BX_J) \cong \pt$.
%We use that $BX_J \longrightarrow BX_I$ is a monomorphism to apply the Sullivan conjecture.
Hence,
we have
a commuting square of weak equivalences
\beq
{
 \xymatrix{
\hcIS \map(\ast, BX_I) \ar[r] \ar[d] & \hcIS \map(BP_n, BX_I)_0 \ar[d]^{(\ref{BLOJAMSeq})}   \\
\map(\ast, BX)\ar[r]                 & \map(BP_n, BX)_0
}
\label{identification}
}
\eeq
so that (\ref{nullcentralizersBP}) is a weak equivalence for all $n$ as desired.
\qed

Given Theorem \ref{nullBX}, the proof of \cite[2.11]{JMratiso} is
directly applicable.  

 \begin{thm}
For $BX$ as in (\ref{BXdefined}) and any homotopy fiber sequence of $p$-compact groups
 \beq
 B\overline{Y} \rightarrow B\widetilde{Y}\stackrel{f}{\rightarrow} BY
  \label{ses}
\eeq
 the natural map induced by $f$ factors as
\beq
{\xymatrix{
& \coprod_{\phi|{B\overline{Y}}\simeq 0} \map(B\widetilde{Y}, BX)_{\phi} \ar[d]  \\
\map(BY, BX) \ar[ur]^\sim \ar[r]^{f^\ast}&  \map(B\widetilde{Y}, BX)
\label{expectionBX}}
}
\eeq
\nod so that $\map(BY, BX)$ is identified, up to homotopy, as the space of maps  $B\widetilde{Y}\stackrel{\phi}{\rightarrow} BY$
 that restrict to nullhomtopic maps from $B\overline{Y}$.
\label{structurethmquo}
\end{thm}
\nod{\bf Proof:} By Theorem  \ref{nullBX},
the natural map
  \beq
\map(\ast, BX) \longrightarrow \map(B\overline{Y}, BX)_0
\label{nullcen}
\eeq
is a weak equivalence.
As in \cite[\S 2]{JMratiso}, we may identify ${f^\ast}$ in (\ref{expectionBX}) with
 the natural map
\beq
\map( BY , \map( \ast , BX)) {\longrightarrow} \map( BY, \map(B\overline{Y}, BX)_0).
\nonumber
\eeq
\qed

With Theorem \ref{structurethmquo} in hand, it is now straightforward to
 adapt %the arguments of
  \cite[\S 6]{JMratiso}.  As first step toward a characterization of nullhomotopic maps from $p$-compact groups $BY$ into $BX$, we consider maps from $p$-toral groups $BP$.

\begin{lem}
Let $BX$ be as in (\ref{BXdefined}) and $BP$ be a $p$-toral group.
A map  $BP \stackrel{f}{\rightarrow} BX$ is nullhomotopic
if and only if for all maps $B\ZZ/p^n\ZZ \stackrel{e}{\rightarrow} BP$
the composition $fe$  is nullhomotopic.
%
%Theorem \ref{nullelememts} holds under the additional hypothesis that $BY$ is
%$p$-toral.
\label{ptoralBXcenter}
\end{lem}
\nod{\bf Proof:} Fix an approximation $\hcl_{n\rightarrow \ff} BP_n \rightarrow BP$.
  For each $BP_n \stackrel{g}{\rightarrow} BX$
 and any $B\ZZ/p^k\ZZ \stackrel{e_n}{\longrightarrow} BP_n$ there is a
homotopy commutative diagram
  \beq
 {\xymatrix{
                        &                                            & BX_I \ar[d]  \\
 B\ZZ/p^k\ZZ \ar[r]^{e_n} & BP_n \ar@{-->}[ur]^{\overline{g}} \ar[r]^g & BX
}
\label{BPn}
\nonumber
}
\eeq
 \nod for some $I \in \cS$ by the identification (\ref{BLOJAMSeq}) for $\pi=P_n$.
  Assuming that
all compositions $ge_n$ are null, then
any $\overline{g}$ is null
by applying \cite[Lemma 6.2]{JMratiso} to $BX_I$.

 Moreover, the identification
(\ref{identification}) implies that
\begin{itemize}
\item  $g$ is null if and only $\overline{g}$ is null and
\item $BP \stackrel{f}{\rightarrow} BX$ is null if and only if $f_n:=f|_{BP_n}$ is null for all $n\in \NN$
\end{itemize}
\nod as the uniqueness obstructions to extending the collection $\{f_n\}_{n \in \NN}$ to  $BP$
lie in $\lim^1(\pi_1(\map(BP_n, BX)_{f_n})=0$.
Hence, if $B\ZZ/p^k\ZZ \stackrel{e}{\rightarrow} BP \stackrel{f}{\rightarrow} BX$ is null
for all $B\ZZ/p^k\ZZ \stackrel{e}{\rightarrow} BP$, then
$f$ is null
by taking $g=f|_{BP_n}$ for varying  $n$.
\qed

   \begin{thm}[cf. {\cite[6.1]{JMratiso}}, {\cite[3.3]{FM2}}]
Let $BX$ be as in (\ref{BXdefined}) and $BY$ be a $p$-compact group.
A map  $BY \stackrel{f}{\rightarrow} BX$ is nullhomotopic
if and only if for all maps $B\ZZ/p^n\ZZ \stackrel{e}{\rightarrow} BY$
the composition $fe$  is nullhomotopic.
\label{nullelememts}
\end{thm}
 \nod{\bf Proof:} Only the ``if'' direction requires proof, so we assume
  all compositions
  $
  B\ZZ/p^n\ZZ \stackrel{e}{\rightarrow} BY \stackrel{f}{\rightarrow} BX
  $
   are nullhomotopic.
Following \cite[\S 6]{JMratiso}, it enough to show that the class $\cat{C}$ of $p$-compact groups such that the statement holds
 has  the five properties of Theorem \ref{saturated}.
  The basic observation, which resulted in the present paper,
 is that arguments to be adapted are not too sensitive to the target of $f$ (see (\ref{piBY'}-\ref{BZBY'}) below).
  Clearly, (i) and (ii) hold.

 Consider property (iii).
 By assumption, all compositions
 \beq
 B\ZZ/p^n\ZZ \stackrel{e}{\rightarrow} BY_0 {\rightarrow} BY \stackrel{f}{\rightarrow} BX
  \nonumber
  \eeq
 are nullhomotopic
 and thus $BY_0 {\rightarrow} BY \stackrel{f}{\rightarrow} BX$ is nullhomotopic.
 Setting $\pi:=\pi_1(BY)$,
 Theorem \ref{structurethmquo} implies that each $f$ factors through $B\pi$, up to homotopy.
 Just as in \cite[Proof of 6.3]{JMratiso},
 we have a lift
  \beq
 {\xymatrix{
B\ZZ/p^m\ZZ \ar@{-->}[r]^{\hat{e}}\ar[d]_{mod ~ p^n} & BY \ar[d] \ar[r]^{f}        & BX  \\
 B\ZZ/p^n\ZZ \ar[r]^{e}                        & B\pi \ar[ur]_{\overline{f}} &
}
\label{piBY'}
}
\eeq
 \nod for {\em any} $B\ZZ/p^n\ZZ \stackrel{e}{\rightarrow} BY$ as indicated.
 Thus, checking property (iii) is reduced to checking that $B\pi \in \cat{C}$.

 Considering property (iv) and setting $BZ$ to be the identity component $\map(BY,BY)_{id}$,
 Theorem \ref{structurethmquo} implies that each $f$  factors through $BZ$, up to homotopy.
  Just as in \cite[Proof of 6.4]{JMratiso},
 we have a lift
  \beq
  {
 \xymatrix{
B\ZZ/p^m\ZZ \ar@{-->}[r]^{\hat{e}}\ar[d]_{mod ~ p^n} & BY \ar[d] \ar[r]^{f}     & BX  \\
 B\ZZ/p^n\ZZ \ar[r]^{e}                        & BZ \ar[ur]_{\overline{f}}&
}
\label{BZBY'}
}
\eeq
 \nod for any $B\ZZ/p^n\ZZ \stackrel{e}{\rightarrow} BY$ as indicated. Thus, the check of  property (iv) is reduced to checking that $BZ \in \cat{C}$.

 Recalling that $BZ$ and $B\pi$ are always $p$-toral groups, Lemma \ref{ptoralBXcenter} above verifies properties (iii) and (iv).

  For property (v), use the homology decomposition (\ref{nullcentralizersdecompostion}) for $BY$.
  In particular, $BY$ being simply-connected and centerless implies each
  \beq
  \cd_\Fp(\map(BV, BY)_v) < \cd_\Fp(BY).
  \nonumber
  \eeq
  By induction,
each $\map(BV, BY)_v \stackrel{v}{\rightarrow} BY \stackrel{f}{\rightarrow} BX$ is nullhomotopic.
  The subspace of $\map(BY, BX)$ consisting of all possible $g$
   such that   all compositions $B\ZZ/p^n\ZZ \stackrel{e}{\rightarrow} BY \stackrel{g}{\rightarrow} BX$ are nullhomotopic
   is weakly equivalent to the subspace
     \beq
\hl_{\cat{A}} \map(\map(BV, BY)_v, BX)_0  .
\nonumber
  \eeq
   \nod To see that this space is path-connected (as desired) apply Theorem \ref{nullBX}.
In particular, $\hl_{\cat{A}} BX$ is weakly equivalent to $BX$ since the $\Fp$-localization of $|\cat{A}^{op}|$ is contractible
and $BX$ is $\Fp$-local.
  \qed

  Finally, we apply Theorem \ref{nullelememts} as in \cite[\S 6]{JMratiso}.

    \begin{thm}[cf. {\cite[6.6, 6.7]{JMratiso}}]
Let $BX$ be as in (\ref{BXdefined})
 and $BY$ be a $p$-compact group with maximal torus $BT \stackrel{\iota}{\rightarrow} BY$
and maximal torus normalizer $BN {\rightarrow} BY$.
For any map $BY \stackrel{f}{\rightarrow} BX$ $f|_{BN}$ is nullhomotopic if and only if $f$ is nullhomotopic.
If---in addition---$BY$ is simply-connected, then
$f|_{BT}$ is nullhomotopic if and only if $f$ is nullhomotopic.
\label{nullconnected}
\end{thm}
\nod{\bf Proof:} See the proofs of Corollaries 6.6 and 6.7 in \cite{JMratiso}. \qed

\begin{remk} As pointed out to the author by Castellana, it is also possible to adapt the arguments of Notbohm \cite{notker}
or to use
the subgroup decomposition for a $p$-compact group \cite{pcptdecomp}---or even the subgroup decomposition for a $p$-local compact group \cite[4.6]{plocalcpt}---to show Theorem \ref{nullBX}.  In these and our approaches, the natural  equivalence
$ \map(B\ZZ/p\ZZ, BX)_0 \simeq BX$ induced by evaluation---obtained via (\ref{BLOJAMSeq})---and the fact that $BX$ is $p$-complete and hence $\FF_p$-local---obtained via Seifert--van Kampen theory---are essential for beginning the inductive process.
%For our specific applications and the proof of Theorem \ref{nullconnected}
%In particular,
\qed
\label{notremark}
 \end{remk}

\section{Null maps between homotopy Kac--Moody groups% and applications
}
\label{sec:nullkacmoody}

Let us now detail some more specific applications. In addition to  $\BKp$,
the results of the previous section apply to the $p$-completed classifying spaces of
\begin{itemize}
\item parabolic subgroups $K_J$  (cf. \cite{classKM}) with $J$ not necessarily spherical  or
\item central quotients (cf. \cite{rank2mv}) and further subgroups of the above % (see \cite{maps}) 
and
\item representatives of $\mathcal{S}^\ast$ defined in \cite{rank2mv} as well as
\item $p$-completed Davis--Januszkiewicz spaces \cites{DJI, DJII}.
\end{itemize}
\nod That is, all these examples can be constructed as in (\ref{BXdefined}) with the $p$-completed classifying space examples
following (\ref{hocoBKlocal}),  representatives of $\mathcal{S}^\ast$ constructed as homotopy pushouts of $p$-completed compact Lie group classifying spaces, and
Davis--Januszkiewicz spaces constructed as a finite poset diagram of torus classifying spaces with an initial object corresponding to the base point.

The next result indicates that the existence of a natural candidate for a maximal torus further simplifies the characterization of
nullhomotopic maps.
 \begin{thm}
For $BX$ and $BX'$ as in (\ref{BXdefined}),
the natural \beq
\map(\ast, BX') \longrightarrow \map(BX, BX')_0
\nonumber
\eeq
\nod is a weak equivalence and
 $BX \stackrel{f}{\rightarrow} BX'$  is nullhomotopic if and only if for all maps $B\ZZ/p^n\ZZ \stackrel{e}{\rightarrow} BX$
the composition $fe$  is nullhomotopic.
If---in addition---the diagram $D$ defining $BX$ has commuting squares
 \beq
 {\xymatrix{
 BT \ar@{=}[r] \ar[d]^{\iota_I} & BT \ar[d]^{\iota_J}   \\
 BX_I \ar[r]^{D(I\le J)} & BX_J
} \label{BXtori}
}
\eeq
\nod for all $I\le J \in \cat{S}$ where $BT \stackrel{\iota_I}{\rightarrow}BX_I$ are maximal tori,
then we have a natural $BT \stackrel{\iota}{\rightarrow} BX$ and
  $BX \stackrel{f}{\rightarrow} BX'$ is nullhomotopic
  if and only if $f|_{BT}$ is nullhomotopic.
\label{nullBX'connected}
\end{thm}
\nod{\bf Proof:} To obtain the natural weak equivalence, we note that Theorem \ref{nullBX} and the contractibility of $|\cat{S}|$ imply that we have
a commuting square of weak equivalences
\beq
{
 \xymatrix{
\map(\ast, BX') \ar[rr] \ar[d] & & \map(BX, BX')_0 \ar[d]^{j_I^\ast}   \\
\map(|\cat{S}|, BX')\ar[rr]& & \hlIS \map(BX_I, BX')_0 .
}
\label{BXtorinull}
}
\eeq

For the first characterization of null maps, only the ``if'' direction requires proof and we assume that
all compositions
 \beq
 B\ZZ/p^n\ZZ \stackrel{e}{\rightarrow} BX_I \stackrel{j_I}{\rightarrow} BX \stackrel{f}{\rightarrow} BX'
  \nonumber
  \eeq
\nod where $e$ is arbitrary and $j_I$ is canonical are null. Theorem \ref{nullelememts} implies that each $fj_I$ is null so that (\ref{BXtorinull}) implies that $f$ is null.

For the second characterization of null maps, only the ``if'' direction requires proof and we assume that $f|_{BT}$ is nullhomotopic.
Recalling that all $BX_I$ are assumed to be simply-connected, Theorem \ref{nullconnected} implies that
each $f|_{BX_I}$ is null so that (\ref{BXtorinull}) implies that $f$ is null.
 \qed

In rough terms, requiring (\ref{BXtori}) means that all $BX_I$ for $I\in \cat{S}$ have a common torus $BT$ which can be added to the diagram $D$ defining $BX$
as an initial space and---since $|\cat{S}|$ is contractible---the homotopy colimit of this extended diagram is still $BX$.

 \nod{\bf Proof of Theorem \ref{intronullconnectedBKp}:} Any $p$-completed Kac--Moody group classifying space can be constructed
  as in (\ref{hocoBKlocal}) so that (\ref{BXtori}) holds for $BT=(BK_\emptyset)\pcp$. \qed

The first three types of examples mentioned at the start of the section have (\ref{BXtori}).
For instance, the pushouts constructing $\mathcal{S}^\ast$ have a rank two torus as their initial spaces
and the two associated  are maps maximal tori. Of course, Davis--Januszkiewicz spaces that are not already the classifying space of a torus do not have (\ref{BXtori}).
We note for completeness that Anjos and Granja \cite{AnjosGranja} constructed a natural space % the space symplectomorphisms of $S^2\times S^2$
as a homotopy pushout of connected Lie group classfying spaces
that does not satisfy (\ref{BXtori}).

 \begin{remk}  A rational analog of Theorem \ref{intronullconnectedBKp} is not possible.
That is, there exist $B\Y \stackrel{f}{\rightarrow} BK^{{\wedge}}_\QQ$ that are not null while $f|_{BT}$ is null.
If $\Y$ is compact connected Lie, then $f|_{BT}$ is null implies that $f$ is null by obstruction theory.  However, the uniqueness obstructions
to extending $0 \simeq f|_{B\Y_I}$ for $I \in \cat{S}$ to $B\Y$ do not vanish in general.  See (\ref{loopsK3}) and Section \ref{sec:higherlimits} below.
Moreover, $\map(B\Y, BK^{{\wedge}}_\QQ)_0$ is not always simply-connected by (\ref{loopsK}).
\qed
 \end{remk}

\section{Applying the arithmetic fiber square}
\label{sec:ratnullkacmoody}

We now begin our applications of Theorem \ref{intronullconnectedBKp} to integral maps $B\Y \stackrel{f}{\rightarrow} BK$ where $\Y$ is a connected compact Lie or
 Kac--Moody group with maximal torus ${\iota}\co T {\rightarrow} \Y$ and $K$ is a Kac--Moody group.
 Because $BK$ is a simply-connected $CW$-complex \cite{nitutkm}, it is given as a homotopy pullback
\beq
{\xymatrix{
BK  \ar[d]^{(-)^{{\wedge}}_\QQ} \ar[rr]^{\prod (-)^{{\wedge}}_p}& & \prod \BKp \ar[d]^{(-)^{{\wedge}}_\QQ}  \\
BK^{{\wedge}}_\QQ  \ar[rr]^{[\prod (-)^{{\wedge}}_p]^{{\wedge}}_\QQ}             & & \prod (\BKp)^{{\wedge}}_\QQ
}
\label{fibresquare}
}
\eeq
\nod known as the arithmetic fiber square \cite{arth} where ${(-)^{{\wedge}}_\QQ}$ denotes rationalization.
By applying $\map(B\Y, -)$ to the above diagram, we have a homotopy pullback of mapping spaces.
 In this section, we obtain
 explicit obstructions to extending a nullhomotopic map $BT \stackrel{0}{\rightarrow} BK$
 to $B\Y \stackrel{f}{\rightarrow} BK$ uniquely.  In other words, we give a description of the set
 homotopy classes of
 $B\Y \stackrel{f}{\rightarrow} BK$ whose restriction to $BT$ is null.
Broading speaking, we now follow the proofs of {\cite[Theorems 3.1 and 3.11(ii)]{JMOself}}.

It turns out that---as in \cite{BKtoBK, JMOself}---$f|_{BT}$ is null implies that ${f}^{{\wedge}}_\QQ$ is null.

   \begin{lem}[cf. {\cite[Proof of Theorem 3.1]{JMOself}}]
The canonical  map
\beq
[B\Y,  BK^{{\wedge}}_\QQ] {\rightarrow}  [B\Y,  \prod (\BKp)^{{\wedge}}_\QQ]
\label{inj}
\eeq
 \nod is an injection.
\label{nullatprimeinj}
\end{lem}
\nod
 Though we are most concerned with the null fiber of (\ref{inj}), all fibers are of general interest.
In particular, if $f\pcp$ is null for all primes $p$, then
\beq
B\Y \stackrel{f}{\rightarrow} BK {\rightarrow} \prod (\BKp)^{{\wedge}}_\QQ
\nonumber
\eeq
\nod is null so that  Lemma   \ref{nullatprimeinj}
implies $f^{{\wedge}}_\QQ$ is nullhomotopic. %---without restriction on $K$ or $\Y$.
Recall that Lemma    \ref{nullatprimeinj} holds whenever  $BK^{{\wedge}}_\QQ$ is a product of $K(\QQ,n)$'s---as in  \cite{BKtoBK, JMOself}.
 In Section \ref{sec:higherlimits} below, we see that $BK^{{\wedge}}_\QQ$ is very often not such a product.  A proof of  Lemma   \ref{nullatprimeinj} is given in Section \ref{sec:vanish}.

Considering Theorem \ref{intronullconnectedBKp} and the preceding discussion, we will now restrict our attention to maps $B\Y \stackrel{f}{\rightarrow} BK$
whose $p$-completions $f^{\wedge}_p$ and rationalization $f^{{\wedge}}_\QQ$ are all nullhomotopic. % with $K$ indecomposable.
 Define $\Nyk \sset [B\Y, BK]$ to be the corresponding set of
homotopy classes.
We can now summarize our discussion in the following technical lemma.

   \begin{lem}
The set $\Nyk$ is in natural bijection with
\beq
\Vyk:= \pi_1(\map(B\Y,  BK^{{\wedge}}_\QQ) , 0) \setminus \pi_1(\map(B\Y,  \prod (\BKp)^{{\wedge}}_\QQ) , 0)
\eeq
\nod so that
 the class of nullhomotopic maps corresponds to $[1] \in \Vyk$ and
any map $B\Y' \stackrel{g}{\rightarrow} B\Y$ induces $\Nyk  \stackrel{g^{\Y}_{\Y'}}{\rightarrow}  \Nyprimek$.  %
\label{nullatprime}
\end{lem}
 \nod{\bf Proof:}
 The homotopy Bousfield--Kan spectral sequence \cite[XI. \S 7]{BK} implies that $\Nyk$ is naturally isomorphic to  $\lim^1\pi_1(\map(B\Y,  P)_0)$ where $P$ is the homotopy pullback diagram (\ref{fibresquare}) for $BK$.
 In particular, $\Nyk$ is naturally isomorphic to  $\Vyk$ since $\pi_1(\map(B\Y,  \prod \BKp ), 0)$ vanishes by Theorem \ref{intronullconnectedBKp}.
\qed

 \nod %We note that
In the more restricted contexts of \cite{BKtoBK, JMOself}, all $\Vyk$ are trivial.
Here we identify
 $\Vyk$ as an (often nontrivial) abelian group so that  any $g^{\Y}_{\Y'}$  corresponds to a group homomorphism. Hence, the kernel of $\Vtk \rightarrow \Vyk$ induced
 by the standard torus
 corresponds to elements of $\Nyk$ which restrict to zero in
$[BT, BK]$. %

First note
 \beq
\pi_1(\map(B\Y,  BK^{{\wedge}}_\QQ) , 0)  &\cong& [B\Y ,\Omega(BK^{{\wedge}}_\QQ) ]
\nonumber \\
                                                         &\cong& [B\Y ,  K^{{\wedge}}_\QQ ]
\nonumber \\
                                                          &\cong& \prod_{j \in J} H^{n_j}( B\Y, \QQ )
%\nonumber
\label{loopsK}
\eeq
\nod where the last identification follows from the fact that $K$ is an $H$-space and a finite type $CW$-complex  (see, e.g., %\cite[3.C]{Hatcher} or
\cite[I: \S2-1]{kane} and \cite[12.2]{nitutkm}).
That is, $H^{\ast}(K, \QQ )$ is free and $K^{{\wedge}}_\QQ$ is a countable product of Eilenberg--MacLane spaces (see also Proposition  \ref{prop:prod} below). Setting $\hat{\QQ}:=(\prod \ZZ\pcp) \otimes \QQ $ we have
 \beq
\pi_1(\map(B\Y,  \prod (\BKp)^{{\wedge}}_\QQ) , 0)  \cong \prod_{j \in J} H^{n_j}( B\Y, \hat{\QQ} )
\nonumber
%\nonumber
\label{loopsK2}
\eeq
\nod as in (\ref{loopsK}).  In conclusion,
 \beq
  \Vyk \cong \prod_{j \in J} H^{n_j}( B\Y, \hat{\QQ}/{\QQ} )
%\nonumber
\label{loopsK3}
\eeq
\nod and we can identify the kernel of $\Vtk \rightarrow \Vyk$.

 \begin{thm}
If $BT \stackrel{B\iota}{\rightarrow} B\Y \stackrel{f}{\rightarrow} BK$ is nullhomotopic, then the group of obstructions
 to $f$ being null is isomorphic to
\beq
\prod_{j \in J_{\rm even}}  \ker(H^{n_j}( B\iota^\ast ,  \hat{\QQ}/\QQ))  \times \prod_{j \in J_{\rm odd}}  H^{n_j}( B\Y ,  \hat{\QQ}/\QQ )
\label{obstructions}
\eeq
 \nod where $J_{\{ \rm odd, even\}}$ is the set of $\{odd,even\}$  $n_j$ (with multiplicity)  from (\ref{loopsK})  and $|J_{\rm odd}|$ is at most the rank of $K$. In particular,
 the set of homotopy classes of such $f$ is isomorphic to the above abelian group.
\label{nullconnectedcorr}
\end{thm}
 \nod{\bf Proof:} Noting that
$H^{\ast}( BT , \hat{\QQ}/\QQ )$
is concentrated in even degrees, it remains only to show that $|J_{\rm odd}|$ is at most the rank of $K$ by the preceding discussion.
Recall the split short exact sequence of Hopf algebras \cites{nitutkm}
 \beq
1 \rightarrow H^{\ast}(K/T_K, \QQ )\otimes_S \QQ \rightarrow  H^{\ast}(K, \QQ ) \rightarrow \Lambda(x_{n_1}, \ldots, x_{n_k}) \rightarrow 1
\label{Hopf}
  \eeq
\nod where $T_K$ denotes the standard torus of $K$, the right algebra is exterior on odd degree generators, $k$ is bounded by the rank of $T_K$,  and the left algebra is concentrated in even degrees (and hence polynomial).
 \qed

 Theorem \ref{nullconnectedLieintro} follows easily from Theorem \ref{nullconnectedcorr}.

  \nod{\bf Proof of Theorem \ref{nullconnectedLieintro}:}
   If $\Y=G$ is a connected compact Lie group with Weyl group $W$ and rank $r$, then   $H^{\ast}(  BG ,   \hat{\QQ}/\QQ )=\hat{\QQ}/\QQ[t_{1}, \ldots, t_{r}]^W$ is concentrated in even degrees
   and $H^{\ast}( B\iota^\ast ,  \hat{\QQ}/\QQ)$ corresponds to the inclusion of $W$--invariants. Thus, the uniqueness obstructions (\ref{obstructions}) vanish as desired. \qed

   As detailed in Section \ref{sec:vanish} below,  $\ker(H^{\ast}( B\iota^\ast ,  \hat{\QQ}/\QQ))$ is the subring of nilpotent elements of $H^{\ast}( B\Y ,  \hat{\QQ}/\QQ)$.
    In this way, Theorem \ref{nullconnectedcorr} describes how obstructions result from the existence nilpotent elements in certain degrees depending on $K$.

 \begin{exmp}
To see that non-trivial obstructions can occur, let $K(A)$ be the   affine Kac--Moody group
associated to the generalized Cartan matrix
  \beq
 A=\left[\begin{array}{ccc} 2 & -1 & -1 \\ -1 & 2 & -1 \\ -1 & -1 &  2 \end{array}\right] .
\nonumber
\label{gcm3}
\eeq
\nod Here it is known that
$
% \beq
 H^{\ast}( K(A), \QQ ) = \Lambda(c_1, y_3, y_5) \otimes \QQ[x_4]
%\nonumber
%\eeq
%\nod
$ as an algebra where subscripts denote degree. %
Using the methods developed in the next two sections, we can calculate
\beq
V^{K(A)}_{K(A)} &\cong& \prod_{i \in \{1, 3, 4, 5\}}H^{i}( BK(A), \hat{\QQ}/\QQ) \nonumber \\
&\cong& H^{4}( BK(A), \hat{\QQ}/\QQ) \times H^{5}( BK(A), \hat{\QQ}/\QQ) \cong (\hat{\QQ}/\QQ)^3 ; \nonumber
\eeq
\nod see Example \ref{exmp:VKA}. This implies
\beq
\ker\left(V^{T}_{K(A)} \rightarrow V^{K(A)}_{K(A)} \right) \cong H^{5}( BK(A), \hat{\QQ}/\QQ) \cong \hat{\QQ}/\QQ. \nonumber
\eeq
\nod Hence, each non-zero $x \in \hat{\QQ}/\QQ$ is represented by a non-trivial $f_x\co BK(A) {\rightarrow} BK(A)$
such that ${f_x}B\iota$ is null.

 We can also construct
a non-empty set of
$f_x\co BK(A) {\rightarrow} BG$ with these properties such that $G$ is compact connected Lie.   That is, $\ker\left(V^{T}_{G} \rightarrow V^{K(A)}_{G} \right)$ is nontrivial
whenever $H^{\ast}(  BG ,\QQ )$ has a polynomial generator in degree 6. Consider, for example, $G=SU(n)$ for $n \ge 3$.
 \qed
\label{exmp:KA}
 \end{exmp}

\section{Vanishing results and the proof of Lemma \ref{nullatprimeinj}}
\label{sec:vanish}

This section gives vanishing results for the left or right terms of (\ref{obstructions}) in Theorem \ref{nullconnectedcorr}
and a proof Lemma \ref{nullatprimeinj}.
%third reviewer comment
 To these ends, we will maintain the notations of the previous section. In particular, $B\Y \stackrel{f}{\rightarrow} BK$
denotes a map between classifying spaces where $K$ is a Kac--Moody group and $\Y$ is a connected compact Lie or
 Kac--Moody group.
Placing conditions on the source $B\Y$ will guarantee that the left term of (\ref{obstructions}) vanishes
while placing conditions on the target $BK$ will guarantee that right term vanishes
so that---under these conditions---$f$ is nullhomotopic if and only if $f|_{BT}$ is null.

  Our line of attack uses the homotopy colimit presentation for $B\Y$.
  Specifically, our discussion here centers around computing
  $[B\Y, BK^{{\wedge}}_\QQ ]$ in terms of $\map(B\Y_I, BK^{{\wedge}}_\QQ)$
  and
$H^\ast(B\Y,  \QQ)$ in terms of $H^\ast(B\Y_I,  \QQ)$ via well-known spectral sequences.
We note that the rational cohomology of the classifying space of a Kac--Moody group is not well-studied %
and refer the reader to
\cite{nitutkm} and Section \ref{sec:higherlimits} for examples.

Let us fix a $BT \stackrel{f}{\rightarrow} BK^{{\wedge}}_\QQ$ and consider the problem of extending $f$ to $B\Y$.
For computation, we appeal to the  homotopy Bousfield-Kan spectral sequence \cite[XI. \S 7]{BK} associated to the functor $\map(B\Y_I, Z)$. It 
has second page given by
\beq
E_2^{i, k} = \lim^i_{\cat{S}^{\rm op}} \pi_k(\map(D_\Y, Z)_{\overline{f}})
\label{homotopy}
\eeq
\nod where we take $Z= BK^{{\wedge}}_\QQ, (\prod \BKp)^{{\wedge}}_\QQ $ and $\overline{f}$ denotes the set of
homotopy classes that extend $f$.
Its related cohomological  spectral sequence \cite[XII. 5.8]{BK} for computing $H^\ast(B\Y,  L)$
has  second page given by
\beq
F_2^{i, k} = \lim^i_{\cat{S}^{\rm op}} H^{k} (B\Y,  L)
\label{cohomology}
\eeq
\nod where we take $L= \QQ ,\hat{\QQ}$. We will use the notation $E_n^{i, k}$ and $F_n^{i, k}$, respectively, for these spectral sequences. The proof of Lemma \ref{nullatprimeinj} depends on comparing them.

In fact, $F_n^{i, k}$ has already been studied in some detail in \cite{nituadams}. % which
%leads directly to following.
 \begin{thm}[cf. {\cite[2.3, 2.5]{nituadams}}]
 For any Kac--Moody group $\Y$ with Weyl group $W$ and $L= \QQ ,\hat{\QQ}$,
 $F_n^{\ast, \ast}$ collapses at $F_2^{\ast, \ast}$ and there are canonical isomorphisms
 \beq
 %fourth reviewer comment
F_2^{i, k} \cong H^i(W,H^k( BT, L)).
\label{cohomologytwo}
\eeq
\nod{}
\label{nituadamsthm0}
\end{thm}
\nod{\bf Proof:}  For the collapse of $F_n^{\ast, \ast}$ at $F_2^{\ast, \ast}$, see the proof of \cite[2.3]{nituadams}.  See the proof of \cite[2.5]{nituadams} for the identification (\ref{cohomologytwo}).
\qed

\nod Moveover, by making suitable identifications, the collapse of $F_\ff^{\ast, \ast}$ implies the collapse of
$E_\ff^{\ast, \ast}$.

 \begin{thm} With the preceding notations, we have
$E_2^{\ast, \ast}=E_\ff^{\ast, \ast}$ for $Z= BK^{{\wedge}}_\QQ ,(\prod \BKp)^{{\wedge}}_\QQ $.
\label{nituadamsthm}
\end{thm}
\nod
Before giving the proof of Theorem \ref{nituadamsthm}, let us first provide further motivation  by relating  the left hand term of (\ref{obstructions}) to $F_2^{i, k}$.

 \begin{prop}
 With the preceding notation, the following are equivalent:
 \begin{enumerate}[(i)]
\item  left hand term of (\ref{obstructions}) vanishes
\item  $H^{n_j} (B\Y,  \QQ)$ has no nilpotent elements for all even $n_j$ from (\ref{loopsK})
\item  $F_2^{2i, n_j-2i}=\lim^{2i} H^{n_j-2i}( D_\Y,  \QQ )=0$ for all even $n_j$ and $i \ge 1$
\end{enumerate}
\label{leftprop}
\end{prop}
\nod{\bf Proof:} %
First note that $H^{\ast} (\iota,  \QQ)$ factors through the edge homomorphism
\beq
 F_\ff^{\ast, \ast} \twoheadrightarrow F_\ff^{0, \ast} \cong H^0(W,H^\ast( BT, \QQ)) \cong  \QQ[t_{1}, \ldots, t_{r}]^W
\label{zerocol}
\eeq
\nod where $\iota$ is the maximal torus of rank $r$.
Because $F_2^{i, k}$ carries a bigraded multiplicative structure,  (ii) and (iii) are equivalent noting Theorem \ref{nituadamsthm0}.  Specifically,
elements in the zero column are not nilpotent and elements in positive columns have nilpotent degree bounded by the maximal chain length in $\cat{S}$.   Likewise, (i) and (iii) are equivalent because the kernel of (\ref{zerocol}) is
precisely the nonzero columns. \qed

Since we are currently focused  on vanishing,  we  recall the following  definition.

\begin{defn}
A Kac--Moody group $K$ is $n$-spherical if in its canonical homotopy colimit presentation (\ref{hocoBK})
the elements of the poset of spherical subsets $\cat{S}$ have cardinality at most $n$.
\label{def:nspherical}
\end{defn}

If $K$ is $n$-spherical, then chains in $\cat{S}$ have length at most $n$
 so that  $\lim^{i>n}_{\cat{S}} M=0$ for {\em any} diagram $M:\cat{S}^{op} \rightarrow \cat{Abelian}$ of abelian groups.
For instance, if $K$ is 3-spherical, then  the last condition in Proposition \ref{leftprop} is simply
to $\lim^2H^{n_j-2}( D_\Y,  \QQ )=0$ for all even $n_j$ from (\ref{loopsK}). Of course, this vanishing is immediate if $K$ is 1-spherical.
More explicitly, if  $\lim^{i\ge2}_{{\cat{S}^{\rm op}}} H^{\ast} (\Y,  \QQ)$, then
$F_2$  yields a short exact sequence
 \beq
 \QQ[t_{1}, \ldots, t_{r}]^W = \lim^0 H^\ast(D_\Y,  \QQ)  \rightarrow  H^\ast(B\Y,  \QQ)   \rightarrow  \lim^1 H^\ast(D_\Y,  \QQ) .
\label{sesBYcoho}
%\nonumber
\eeq
\nod As we have seen, this implies that the left hand term of (\ref{obstructions}) vanishes. %(and  the complexity of computing $H^\ast(B\Y,  \QQ)$ is low).

If we can replace $D_\Y$ with a subdiagram over a subposet with a smaller chain length bound,
then this subdiagram provides a more favorable $E_2$.  For this purpose, we make the following definition.
\begin{defn}
A Kac--Moody group $\Y$ is essentially $n$-spherical if the exists a right cofinal full subposet $\cat{S}'\sset \cat{S}$
whose chains have length at most $n$.
\end{defn}
In lieu of a definition of right cofinal (see, e.g. \cite[19.6.1]{Hirschhorn}), we give the two properties relevant to our discussion
and
a practical criterion for checking if $\Y$ is essentially $n$-spherical.  Specifically,  a right cofinal $\cat{S}'\sset \cat{S}$
induces a canonical homotopy equivalence (see, e.g. \cite[19.6.17]{Hirschhorn})
 \beq
 \hcl_{I\in \cat{S}'} B\Y_I \stackrel{\thicksim}{\longrightarrow} \hcl_{I\in \cat{S}} B\Y_I \stackrel{\thicksim}{\longrightarrow} BK
\label{hocolimcompare}
%\nonumber
\eeq
\nod so that Theorem \ref{nituadamsthm0} implies that the canonical comparison maps %
%as well as canonical isomorphisms (cf. \cite[pp. 153--157]{GZ})
 \beq
 \lim^i_{(\cat{S}')^{\rm op}} H^\ast(D_\Y,  \QQ)  \stackrel{\sim}{\longrightarrow}   \lim^i_{\cat{S}^{\rm op}} H^\ast(D_\Y,  \QQ)
\label{limcompare}
%\nonumber
\eeq
\nod
are isomorphisms for all $i \ge 0$ (see also, e.g., \cite[1.5]{cofinal}). The particularities of our situation yield the following. % proposition gives an easily checkable criterion.
 \begin{prop}
For a Kac--Moody group $K$ with spherical subsets $\cat{S}$,
there exists a unique minimal, right cofinal, full subposet ${\overline{\cat{S}}}\sset \cat{S}$
whose elements are intersections of maximal elements of $\cat{S}$.
In particular, $K$ is essentially $n$-spherical if and only if
all sets of $n$ distinct maximal elements of the poset of spherical subsets have a common intersection.
\label{rightcofinalprop}
\end{prop}
\nod{\bf Proof:} For any finite poset, the existence of a unique minimal, right cofinal, full subposet is guaranteed by \cite[Theorem 4.3]{MoGel}. For $\cat{S}$, it is convenient identify this subposet by considering
covers of categories \cite[2.5]{PhDart}.

 To see that ${\overline{\cat{S}}}\sset \cat{S}$ is right cofinal, note that
the maximal elements provide a canonical cover of $\cat{S}$ and $\cat{S}$ has greatest lower bounds.  That is, the collection of $\cat{S} \downarrow I_m \sset \cat{S}$ defined as the
full subposet of elements under a fixed maximal $I_m \in \cat{S}$ is a cover and ${\overline{\cat{S}}}\sset \cat{S}$ factors as a composition of cofinal functors
 ${\overline{\cat{S}}} \rightarrow {\overline{\cat{S}}} \ltimes \cat{S} \downarrow I \rightarrow \cat{S}$ \cite[Proposition  2.5.3]{PhDart}.

 To see that  $\overline{\cat{S}}$ is minimal,
we must show that $I \notin \overline{\cat{S}}$ implies that the full subposet of elements strictly greater than $I$ denoted $I \downarrow \downarrow\cat{S} $ has a contractible realization \cite[Theorem 4.3]{MoGel}.
Since $I \notin \overline{\cat{S}}$, the maximal elements of $I \downarrow  \downarrow\cat{S} $ constitute a non-empty subset of the maximal elements of $\cat{S}$ and provide a cover of
$I \downarrow \downarrow\cat{S} $ as above.  The analogous composition of cofinal functors now shows that the realization of $I \downarrow \downarrow\cat{S} $ is homotopic to a $(k-1)$--simplex where
$k$ is the number of maximal elements in $I \downarrow  \downarrow\cat{S} $.

%If we take $\cat{S}$ to be the poset of spherical subsets associated to a Kac--Moody group $K$, then
Clearly, the maximal chain length $l$ of $\overline{\cat{S}}$ is the minimal $n$ such that  $K$ is essentially $n$-spherical.
The initial object of $\overline{\cat{S}}$ is
necessarily the intersection of all maximal elements of $\cat{S}$ and must coincide with the intersection of
all sets of $l$ distinct all maximal elements of $\cat{S}$.
\qed

Let us now consider the right hand term of (\ref{obstructions}).  We already have seen that $J_{\rm odd}$ is finite
and when $K$ is indefinite indecomposable, we have further control on $J_{\rm odd}$.
 \begin{prop}
If $K$ is indefinite indecomposable, then
the right hand term of (\ref{obstructions}) vanishes---without further restrictions on $\Y$.
\label{rightprop}
\end{prop}
 \nod{\bf Proof:} For indefinite, indecomposable $K$,  %neither affine nor Lie,
 odd $n_j$ equal one or three by \cite[Appendix]{chow}, \cite{ratcohoWeyl} (see also Appendix \ref{append}).
Since $B\Y$ is simply-connected, it remains to show
 that $H^{3}( B\Y , \QQ )=0$.
By applying the the Rothenberg--Steenrod spectral sequence  to $H^{\ast}(\Y , \QQ )=0$ (see, e.g., \cite[II: \S6-2 Theorem B (ii)]{kane}), this vanishing is easily verified by recalling (\ref{Hopf}) and
 checking that $H^{2}(\Y/T, \QQ )\otimes_S \QQ=0$.
 This latter vanishing follows, for example, from the fact that degree two elements in $S^2$ generate
  %fifth reviewer comment
$H^{2}(\Y/T, \QQ )$ \cite[\S
6]{nitutkm}. % (cf. ).% or $H^{3}(B\Y, \QQ )\cong \lim^1 H^{2}(D, \QQ )=0$ directly.
\qed

We can now summarize our discussion with the following theorem.

 \begin{thm}
 If $K$ is indefinite, indecomposable, then
 the following are equivalent:
 \begin{enumerate}[(i)]
\item  $B\Y \stackrel{f}{\rightarrow} BK$ is null if and only if $f|_{BT}$ is null
\item  $H^{n_j} (B\Y,  \QQ)$ has no nilpotent elements for all even $n_j$ from (\ref{loopsK})
\item   $F_2^{2i, n_j-2i}=\lim^{2i} H^{n_j-2i}( D_\Y,  \QQ )=0$ for all even $n_j$ and $i \ge 1$
\end{enumerate}
In particular, these conditions are equivalent to $\lim^2H^{n_j-2}( D_\Y,  \QQ )=0$
when $\Y$ is essentially 3-spherical and always
hold when $\Y$ is essentially 1-spherical.
\label{selfmaps}
\end{thm}
\nod   The proof of
Theorem  \ref{nullconnectedintro} has now been reduced to a nilpotency condition on rational cohomology.
Over the next two sections, we will prove the following.

 \begin{thm}
 If a Kac--Moody group $K$ is 2-spherical, then
$\lim^2 H^\ast( D_K,  \QQ )=0$ where $D_K$ is the diagram for its canonical homotopy colimit presentation (\ref{hocoBK}).
\label{lim2vanish}
\end{thm}

\nod Hence, combining the above with Theorem \ref{selfmaps} proves Theorem  \ref{nullconnectedintro}.  More generally, the following is immediate.
  \begin{thm}
Under either of the following two additional hypotheses:
\begin{itemize}
\item  $\Y$ is connected compact Lie
\item $\Y$ is  a 2-spherical Kac--Moody group and $K$ has no affine or Lie direct product factors,
\end{itemize}
${f}\co B\Y {\rightarrow} BK$ is nullhomotopic $\iff$ $f|_{BT}$ is nullhomotopic.
\label{nullconnectedmain}
\end{thm}
\nod To see that this hypothesis on $K$ is necessary, note Example \ref{1spherical} below.

We now turn to the proof of Lemma    \ref{nullatprimeinj} which is  much
 easier if $\Y$ is essentially 1-spherical and motivates our approach.
 \begin{prop}
  %sixth reviewer comment
If  $\Y$ is essentially 1-spherical,
then (\ref{inj}) is an injection.
\end{prop}

\nod{\bf Proof:} As in (\ref{sesBYcoho}), we have the collapse at $E_2$ of the associated homotopy spectral sequence (\ref{homotopy}) and  we obtain a map of short exact sequences
 \beq
 {\xymatrix{
 \lim^1_{\overline{\cat{S}}^{\rm op}} \pi_1(\map(D_\Y,  BK^{{\wedge}}_\QQ)) \ar[r] \ar[d]            & [B\Y,  BK^{{\wedge}}_\QQ]  \ar[r] \ar[d]               & \lim_{\overline{\cat{S}}^{\rm op}} [D_\Y,  BK^{{\wedge}}_\QQ]   \ar[d]      \\
 \lim^1_{\overline{\cat{S}}^{\rm op}} \pi_1 (\map(D_\Y,  \prod (\BKp)^{{\wedge}}_\QQ)) \ar[r]         & [B\Y,  \prod (\BKp)^{{\wedge}}_\QQ]  \ar[r]            & \lim_{\overline{\cat{S}}^{\rm op}} [D_\Y,  \prod (\BKp)^{{\wedge}}_\QQ].
}\label{sesBY}
\nonumber
}
\eeq
\nod In particular, it is sufficient to show that the leftmost vertical map is an injection because the rightmost
vertical map corresponds to an inclusion of $W$-fixed points
\beq
\QQ[t_{1}, \ldots, t_{r}]^W \hookrightarrow \hat{\QQ}[t_{1}, \ldots, t_{r}]^W.
\nonumber
\eeq
By (\ref{loopsK})  and (\ref{sesBYcoho}), %computing the cohomology of $B\Y$ in terms of $B\Y_I$,
this map of $\lim^1$ terms corresponds to a change in coefficients
\beq
\prod_{j \in J_{\rm even}}  H^{n_j+1}( B\Y ,  \QQ \hookrightarrow \hat{\QQ} ).
\nonumber
\eeq
\qed

Using Theorem \ref{nituadamsthm0}, this argument can be refined to handle the general case.

\nod{\bf Proof of Theorem \ref{nituadamsthm}:}  %For the collapse of $F_n^{\ast, \ast}$ at $F_2^{\ast, \ast}$, see the proof of \cite[2.3]{nituadams}.  See the proof of \cite[2.5]{nituadams} for the identification \ref{cohomologytwo}.
Consider the $E_2^{\ast,0}$. The uniqueness obstructions to extending $f$ to $B\Y_I$ are canonically isomorphic
to $H^i(W_I; \pi_i(\map(BT, Z)_f)$ and vanish.  Hence  $E_2^{i,0}=0$ for $i \ge 1$.

 Now let $F_n^{\ast, \ast}$ be the associated spectral cohomology sequence (\ref{cohomology}) with coefficients in  $L=\QQ, \hat{\QQ}$.
For $n \ge 2$, the differentials $e_n$ on $F_n^{\ast, \ast}$
%seventh reviewer comment
are closely related to the differentials $d_n$ on $E_n^{\ast, \ast}$ by (\ref{loopsK}). That is,  $E_2^{i, k}$
is given by
\beq
\lim^i \pi_k(\map(D_\Y, Z)) \cong \lim^i \prod_{j \in J_{\rm odd}}  H^{n_j-k+1}( D_\Y,  L ) \cong \prod_{j \in J_{\rm odd}} \lim^i   H^{n_j-k+1}( D_\Y,  L )
\nonumber
\eeq
\nod for $k\ge 1$ whereas $F_2^{i, k}$ is given by
$
\lim^i   H^{j}( D_\Y ,  L )
$.
 By definition \cite[XII. 5.8]{BK}, $d_2^{i, k}$  is given as a product of $e_2^{i, n_j-k+1}$ wherever the source or
 target is not in the zero column. Equivalently, we can consider (the shift of) the spectrum associated to $K^{{\wedge}}_L$.  In particular, $d_n^{i, k}=0$ for all $n\ge2$
 by the collapse of $F_n^{\ast, \ast}$ ($k\ge 1$) and the first paragraph of this proof ($k=0$).
 \qed

\nod{\bf Proof of Lemma    \ref{nullatprimeinj}:}
  The map (\ref{inj}) is induced by a map on
$E_2=E_\ff$ (\ref{homotopy}) % of the associated Bousfield-Kan spectral sequences \cite{BK} which abuts to $\pi_i(\map( B\Y, Z)$
associated to $ BK^{{\wedge}}_\QQ \rightarrow (\prod \BKp)^{{\wedge}}_\QQ $. Moreover, (\ref{inj}) is an injection
if and only if the restriction of this map of $E_\ff$-pages to total degree zero is an injection.

The change in coefficients $\QQ \hookrightarrow \hat{\QQ}$ induces  an injection on cohomology. Therefore, the map induced on $F_2=F_\ff$ is an injection.    By the
above identifications,  the corresponding map of the total degree zero parts of $E_\ff$-pages must be an injection.  Thus, (\ref{inj}) is an injection as desired.
\qed

\section{Explicit $\lim^2$ computations and examples}
\label{sec:higherlimits}

 If a Kac--Moody group $K$ is 2-spherical, then we saw in the previous section that---since the maximal chain length in $\cat{S}$ is at most two---Theorem  \ref{nullconnectedintro}
 is equivalent to the vanishing of $\lim^2H^{n_j-2}( D_K,  \QQ )$   for all even $n_j$ from (\ref{loopsK}) (see Theorem \ref{selfmaps}).
In this section, we give a method  to compute $\lim^2_\cat{S} M$ for an arbitrary diagram of free $R$-modules
$M\co\cat{S}^{op} \rightarrow \cat{$R$-mod}$ over some fixed ring $R$.

This method is applied in the next section to prove  Theorem \ref{lim2vanish} and obtain the desired vanishing.
For 2-spherical $K$,  Theorem \ref{lim2vanish} also implies that the short exact sequence (\ref{sesBYcoho}) is
available to compute $H^\ast( BK,  \QQ )$.  The second purpose of this section is to give examples where
(\ref{sesBYcoho}) computes $H^\ast( BK,  \QQ )$ as a ring. Let us now give a 1-spherical example.
%as a concrete illustration of how (\ref{sesBYcoho}) can be applied.
%
 \begin{exmp}
Taking, for example, the  non-symmetrizable, non-singular matrix
  \beq
 A_1=\left[\begin{array}{ccc} 2 & -1 & -1 \\ -7 & 2 & -1 \\ -8 & -9 &  2 \end{array}\right] ,
\nonumber
\label{gcm3}
\eeq
$K:=K(A_1)$ is indefinite, indecomposable, and 1-spherical with rank 3 torus.   These conditions imply
that $\QQ[t_{1}, t_2, t_{3}]^W= \QQ$ (see Table \ref{table:Jodd}) and
 $H^\ast(BK,  \QQ)$ is necessarily concentrated in odd degrees by (\ref{sesBYcoho}).
Moreover, a simple Euler characteristic calculation gives the dimension of each odd degree.  Specifically, it is straightforward to compute the Poincar\'{e} series of each degree
of the standard normalized chain complex whose homology computes $\lim^i_\cat{S} H^\ast(BK_I,  \QQ)$ (see, e.g., \cite[Appendix II.3]{GZ} and \cite[8.3]{Weibel} or \cite[2.3]{CGNDMT}). Because $\lim^0_\cat{S} H^{\ast\ge 1}(BK_I,  \QQ)=0$, the degree zero and one Poincar\'{e} series
give the Poincar\'{e} series for $\lim^1_\cat{S} H^\ast(BK_I,  \QQ)$.  In this way,  the Poincar\'{e} series for $H^{\ast\ge1}(BK,  \QQ)$ is given by:
\beq
t + \frac{2t}{(1-t^2)^3}-\frac{3t}{(1-t^2)^2(1-t^4)}
&=& t + \sum_{n=1}^{\ff} \frac{2n^2-5-3(-1)^{n}}{8}   t^{2n+1} \nonumber\\
&=& 2t^{7} + 3t^{9} + 6t^{11} + 8t^{13} + 12t^{15} + 15t^{17} + \cdots
\nonumber
\eeq
\nod
In particular, $H^\ast(BK,  \QQ)$ is infinitely generated  and  $BK^{{\wedge}}_\QQ$ is not a product of Eilenberg--MacLane spaces.  Indeed, the product of any two odd degree elements has positive, even degree and necessarily vanishes.
This determines $H^\ast(BK,  \QQ)$ as a ring.
  %7th reviewer comment
The careful reader will have already noted the specific form of $A_1$ is not needed for this computation; the stated properties (indefinite, indecomposable, and 1-spherical with rank 3 torus) are sufficient.

We emphasize that with $H^\ast(BK,  \QQ)$ in hand we can read off the set of homotopy classes of
$f_x\co BK {\rightarrow} BG$ such that $G$ is compact connected Lie
and ${f_x}B\iota$ is null from (\ref{obstructions}). For instance, if  $G=SU(n+1)$, then $G\ratcp \simeq K(\QQ, 3) \times \cdots \times K(\QQ, 2n+1)$ and the group of uniqueness obstructions given by
(\ref{obstructions}) is
\beq
H^3(BK,  \hat{\QQ}/\QQ) \times \cdots \times H^{2n+1}(BK,  \hat{\QQ}/\QQ) %
\cong (\hat{\QQ}/\QQ)^{\lfloor \frac{2n^3+3n^2-4{n}}{24} \rfloor + (-1)^{n+1} } \nonumber
\eeq
\nod
That is, up to isomorphism, the sequence of obstruction groups begins
\beq
\left \{ \ker \left(V^{T}_{SU(n+1)} \rightarrow V^{K}_{SU(n+1)} \right) \right \}_{n \ge 1} = 0,0, (\hat{\QQ}/\QQ)^2, (\hat{\QQ}/\QQ)^5, (\hat{\QQ}/\QQ)^{11}, (\hat{\QQ}/\QQ)^{19},  \ldots %(\hat{\QQ}/\QQ)^{31}, \ldots
\nonumber
\eeq
\nod
 Similarly, it is straightforward to exhibit affine Kac--Moody groups $K'$ such that homotopically nontrivial $f_x\co BK(A_1) {\rightarrow} BK'$ with ${f_x}B\iota$ nullhomotopic exist.
 Consider, for example, the (untwisted) affine Kac--Moody group $K'$ obtained by extending $SU(n+1)$ for $n\ge3$ \cite[XIII]{kumar}.
\qed
\label{1spherical}
 \end{exmp}

  As in  Example \ref{1spherical}, our computation of  $\lim^2_{\cat{S}} M$ for some  diagram of free $R$-modules $M\co\overline{\cat{S}}^{op} \rightarrow \cat{$R$-mod}$ begins with
 the standard normalized chain complex associated to $M$.   First, we give a concrete description of $\lim^2_{\cat{S}} M$ in a special case.
Second, we show that $\lim^2_\cat{S} M$ surjects onto $\prod \lim^2_{\cat{S}_k} M$ for a finite collection of subposets $\cat{S}_k \sset \cat{S}$
such that: (1) terms of this product can be  described as in the special case; and (2) the kernel of this surjection is explicit.

Our arguments employ algebraic discrete Morse theory and depend only on the form of $\cat{S}$.   In particular,
they are relevant to essentially 2-spherical Kac--Moody groups as well as situations unrelated to Kac--Moody groups where $\lim^2 M$ is nonzero.  For further details on (algebraic) discrete Morse theory, we give
specific references in \cites{DMT,aDMT,aDMTlong,CGNDMT}.  We are pleased to recommend \cite[Sections 2-3]{CGNDMT} in particular for
  a reasonably self-contained introduction that is formulated cohomologically. Our applications of discrete Morse theory will  be fairly elementary and motivating figures will be given.

The following lemma constitutes the first step. % and we include its proof for completeness.

 \begin{lem}
Define the poset of subsets of $\{ 1, \ldots , n\}$
\beq
\cat{C}_n=\{ \emptyset , \{ 1 \}, \ldots , \{n\}, \{ 1, 2 \}, \ldots , \{n-1, n\}, \{n, 1\} \}
\eeq
\nod ordered by inclusion.  For any $n \ge 3$ and diagram of free $R$-modules, $M\co\cat{C}_n^{op} \rightarrow \cat{$R$-mod}$
\beq
\lim^2_{\cat{C}_n} M\cong M_\emptyset / \langle M_1^\emptyset, \cdots, M_n^\emptyset \rangle
\label{inv2star}
%\nonumber
\eeq
\nod where
$M_I:=M(I)$, $M_i:=M(\{ i \})$, and $M_i^\emptyset :=\im (M(\{ i \} \rightarrow \emptyset)) $.
\label{lim2special}
\end{lem}
\nod{\bf Proof:}
Here the standard normalized chain complex has the form
 \beq
0 \rightarrow C^0 \rightarrow C^1 \rightarrow C^2 \rightarrow 0
\label{chainhstart}
\eeq
\nod  with $C^n = \prod_{ I_0 \varsubsetneq \cdots \varsubsetneq I_n} M(I_0 ) $
and is chain homotopic to
 \beq
0 \rightarrow C^0 / M_\emptyset %\prod_{I \in \cat{C}_n, |I|=2} M_I
  \rightarrow \prod_{ 1 \le i \le n} M_i \times M_i  \rightarrow M_\emptyset \rightarrow 0
\label{chainhend}
\eeq
\nod  so that $H^2 (C^\ast)$ is the cokernel of  $ M_1 \times M_1 \times \cdots  \times M_n \times M_n \rightarrow M_\emptyset$
given by
\beq
(x_1, \overline{x}_1 , \ldots, x_n , \overline{x}_n ) \mapsto x_1^\emptyset + \overline{x}_1^\emptyset + \cdots + x_n^\emptyset + \overline{x}_n^\emptyset
\label{diff}
\eeq
\nod   where  $x_i^\emptyset, \overline{x}_i^\emptyset \in M_{\emptyset}$ are the respective images of $x_i, \overline{x}_i \in M_i$  under  $M(\{ i \} \rightarrow \emptyset)$.

\begin{figure}[!h]
\centering
\begin{tikzpicture}[scale=.8]
\def\initalrotate{0}
\def\boundarysimplex#1#2{
  \draw #1 #2 +(30+\initalrotate:2)    node[anchor=\initalrotate+210]{$3$}
 -- node[rotate=0+\initalrotate] {\midarrow}  +(90+\initalrotate:1) node[anchor=\initalrotate+270]{$13$}
   --node[rotate=180+\initalrotate] {\midarrow} +(150+\initalrotate:2) node[anchor=\initalrotate+330]{$1$}
  -- node[rotate=120+\initalrotate] {\midarrow} +(210+\initalrotate:1) node[anchor=\initalrotate+390]{$12$}
 -- node[rotate=300+\initalrotate] {\midarrow} +(270+\initalrotate:2)  node[anchor=\initalrotate+450]{$2$}
 -- node[rotate=240+\initalrotate]{\midarrow} +(330+\initalrotate:1) node[anchor=\initalrotate+510]{$23$}
 -- node[rotate=60+\initalrotate] {\midarrow} +(390+\initalrotate:2);
}

\def\simplex#1#2{
  \draw #1 #2 +(30+\initalrotate:2)    node[anchor=\initalrotate+202]{$12$}
   --node[anchor=270+\initalrotate] {$1$} +(150+\initalrotate:2) node[anchor=\initalrotate+338]{$13$}
 -- node[anchor=30+\initalrotate] {$2$} +(270+\initalrotate:2)  node[anchor=\initalrotate+450]{$23$}
 -- node[anchor=510+\initalrotate] {$3$} +(390+\initalrotate:2);
}

\def\rightmapstoarrow#1#2{
 \draw [|->] #1 #2 -- +(1,0)  ;
}

\def\simplexarrowa#1#2{
  \draw #1 #2
 -- node[rotate=210+\initalrotate] {\midarrow}  +(30+\initalrotate:2)   ;
}

\def\simplexarrowapairing#1#2{
\draw #1 #2 [->,red, thick] ++(30+\initalrotate:.8) -- +(120+\initalrotate:.2);
}

\def\simplexarrowb#1#2{
  \draw #1 #2
 -- node[rotate=270+\initalrotate] {\midarrow}  +(90+\initalrotate:1) ;
}

\def\simplexarrowbpairing#1#2{
\draw [->,red, thick]  #1 #2 ++(90+\initalrotate:0.35) -- +(180+\initalrotate:0.2) ;
}

\def\dotpairing#1#2{
\draw [->,green, thick]  #1 #2 +(0:0) -- +(30+\initalrotate:0.2) ;
}

\def\dotpairingedges#1#2{
\draw [->,green, thick]  #1 #2 ++(30+\initalrotate:2) -- +(240+\initalrotate:0.2) ;
}

  \def\3dotunlabeledsimplex#1#2{
\draw[fill=black] #1 #2 + (30:2) circle (.6mm)
   +(150:2) circle (.6mm)
 +(270:2)  circle (.6mm) -- cycle;
  }

    \def\6dotunlabeledsimplex#1#2{
\draw[fill=black] #1 #2 + (30:2) circle (.6mm)
  +(90:1) circle (.6mm)
   +(150:2) circle (.6mm)
   +(210:1) circle (.6mm)
 +(270:2)  circle (.6mm)
  +(330:1) circle (.6mm) -- cycle;
  }

    \def\1dotunlabeledsimplex#1#2{
\draw[fill=black] #1 #2 +(0:0) circle (.6mm);
  }

\boundarysimplex{(0,0)}{(-2,0)}

\6dotunlabeledsimplex{(0,0)}{(-2,0)}
\1dotunlabeledsimplex{(0,0)}{(-2,0)}

\rightmapstoarrow{(0,0)}{(0,-0.5)}

  \foreach \a in {0, 120, 240} {
  \def\initalrotate{\a}
  \simplexarrowa{(0,0)}{(-2,0)}
    \simplexarrowb{(0,0)}{(-2,0)}
    \simplexarrowbpairing{(0,0)}{(-2,0)}
    }

      \foreach \a in {240, 120} {
  \def\initalrotate{\a}
    \simplexarrowapairing{(0,0)}{(-2,0)}
    }
        \foreach \a in {0} {
  \def\initalrotate{\a}
   \dotpairing{(0,0)}{(-2,0)}
    }

    \boundarysimplex{(0,0)}{(3,0)}

\6dotunlabeledsimplex{(0,0)}{(3,0)}

  \foreach \a in {0, 120, 240} {
  \def\initalrotate{\a}
\dotpairingedges{(0,0)}{(3,0)}
    }

\rightmapstoarrow{(0,0)}{(5.5,-0.5)}

    \simplex{(0,0)}{(8,0)}

\3dotunlabeledsimplex{(0,0)}{(8,0)}

\end{tikzpicture}

\caption{A picture of our discrete Morse theory pairing for $\cat{C}_3$. }
\label{pic:C3}
\end{figure}
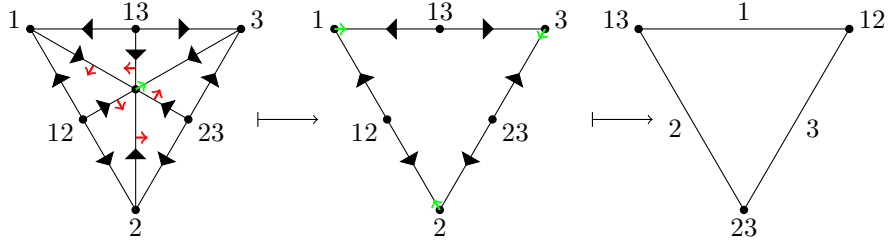

   From the perspective of algebraic discrete Morse theory, this chain homotopy is induced by a single acyclic matching \cite[2.1]{aDMTlong}, \cite[2.4]{CGNDMT} between $M_{\emptyset}$ factors.  In ordinary discrete Morse theory \cites{DMT}, pairs of cells in a (simple) $CW$-complex can be eliminated by elementary collapse and an acyclic matching combines multiple pairings into a single homotopy equivalence. For instance, if we ignore the triangular arrows in Figure \ref{pic:C3}, then the first step depicts an acyclic matching with six pairs of cells corresponding to the remaining six arrows, colored red and green.  The second step in Figure \ref{pic:C3} depicts an acyclic matching with three pairs of cells.  Continuing in this way, we can describe an explicit homotopy from the starting complex to a point.   In what follows, we will collapse (geometrically motivated) pairs of free $R$-modules via acyclic matchings in the sense of algebraic discrete Morse theory.
   
     Now let us explicitly describe the single acyclic matching inducing the chain homotopy from (\ref{chainhstart}) to (\ref{chainhend}).
   The unique $M_{\emptyset}$ factor in $C^0$ is paired with the $M_{\emptyset}$ factor in $C^1$ associated with the morphism $\emptyset \rightarrow n$ and all other $M_\emptyset$ factors in $C^1$ factors are paired with
   $M_{\emptyset}$ factors in $C^2$  as follows: each $\emptyset \rightarrow l$ coordinate for $l <n$ pairs with
  the $\emptyset \rightarrow l+1 \rightarrow \{l, l+1\}$ coordinate, each $\emptyset \rightarrow \{l, l+1\}$ coordinate pairs with the
   $\emptyset \rightarrow l \rightarrow \{l, l+1\}$ coordinate, and the $\emptyset \rightarrow \{1, n\}$
   coordinate pairs with the
   $\emptyset \rightarrow 1 \rightarrow \{1, n\}$ coordinate.  Strictly speaking, this pairing involves choosing a basis for $M_{\emptyset}$ and pairing these
bases identically. It is easy to check that this an acyclic pairing.  The differential (\ref{diff}) is an instance of \cite[2.2]{aDMTlong}, \cite[2.5]{CGNDMT}.
   \qed

   When $n=3$, the acyclic matching just described is
depicted by the first step of Figure \ref{pic:C3} with the six non-triangular arrows representing six pairings of $M_\emptyset$ factors.
 Intuitively, removing these 12 copies of $M_\emptyset$ from the chain complex (\ref{chainhstart}) is an algebraic version of a homotopy equivalence that eliminates 12 cells.

%%%%%%%%%%%%%%%%%%%%%%%%%%%%%%%%%%%%%%%%%%%%%%%%%%%%%%%%%

It is well-known that the chain complex can be further reduced. The second step in Figure \ref{pic:C3} indicates an algebraic discrete Morse theory proof, but we do not spell this out.  Performing only the  first step is better motivation for later computations.
 Let us now give a 2-spherical example where Lemma \ref{lim2special} applies.

 \begin{exmp}
Taking, for example, the  non-symmetrizable, non-singular matrix
  \beq
 A_2=\left[\begin{array}{ccc} 2 & -1 & -1 \\ -1 & 2 & -1 \\ -2 & -3 &  2 \end{array}\right] ,
\nonumber
\label{gcm3}
\eeq
$K:=K(A_2)$ is indefinite, indecomposable, and 2-spherical with rank 3 torus $T:=K_\emptyset$.  % These conditions imply
%that $\QQ[t_{1}, t_2, t_{3}]^W= \QQ$ and
In particular, $\cat{C}_3$ is the set of spherical subsets $\cat{S}$ associated to $K$
so that Lemma \ref{lim2special} gives the identification
\beq
\lim^2_{\cat{S}} H^n( D_{K},  \QQ )= H^n( BT,  \QQ ) / \langle H^n( BK_{\{1\}},  \QQ ), \cdots, H^n( BK_{\{3\}},  \QQ)  \rangle
\nonumber
\eeq
\nod
In order to use the short exact sequence (\ref{sesBYcoho}) as in Example \ref{1spherical},  it remains to show that
the $\QQ$-linear span of $\{ H^n( BK_{\{i\}},  \QQ ) \}_{i=1,2,3}$ is $H^n( BT,  \QQ )$.

Since
the Weyl group of $K$ has an infinite set of reflecting hyperplanes, we can apply
Lemma \ref{lim2reflect}
 with
$s$ generating the Weyl group of $K_1$, $t$ generating the Weyl group of $K_2$ and $u$ generating the Weyl group of $K_3$.  (For reference, a discussion of the Weyl group of a Kac--Moody group
is given in Section \ref{sec:invarianttheory}.)
In particular,
$\{ H^n( BK_{\{i\}},  \QQ ) \}_{i=1,2,3}$ spans $H^n( BT,  \QQ )$ as desired.

So we have (\ref{sesBYcoho}) as in Example \ref{1spherical}.
Hence, $H^{\ast \ge 1}(BK,  \QQ)$ is necessarily concentrated in odd degrees, the cup product of any two odd degree elements is zero, and the rank of each odd degree is determined by
an Euler characteristic calculation.  In this way,  the Poincar\'{e} series for $H^{\ast}(BK,  \QQ)$ is given by:
\beq
1 &+&   t \biggl(1 + \frac{-1}{(1-t^2)^3} + \frac{3}{(1-t^2)^2(1-t^4)} - \frac{1}{(1-t^2)(1-t^4)(1-t^6)} \nonumber \\
 &-&  \frac{1}{(1-t^2)(1-t^4)(1-t^8)}-\frac{1}{(1-t^2)(1-t^4)(1-t^{12})}\biggr) \nonumber \\
  &=& 1 + t^{7} + t^{9} + 3t^{11} + 2t^{13} + 5t^{15} + 7t^{17} + 7t^{19} +\cdots %
\nonumber
\eeq
\nod and only one ring structure on $H^{\ast}(BK,  \QQ)$ is possible. %
Here the specific form of $A_2$ effects the result, but the same method computes $H^\ast(BK,  \QQ)$ as a ring for any $A_2$ with the properties stated  above.
As in Example \ref{1spherical}, $f\co BK \rightarrow BK$  is null only if ${f}|_{BT}$ is null,  but it is straightforward to exhibit  affine Kac--Moody $K'$  or compact connected Lie $G$ such that homotopically nontrivial $g\co BK {\rightarrow} BK'$  and $h\co BK {\rightarrow} BG$ with $g|_{BT}$ and $h|_{BT}$  nullhomotopic exist.
\qed
\label{2spherical}
\label{2sphericalcont}
 \end{exmp}

Let us now consider the form of $\cat{S}$ associated to a 2-spherical Kac--Moody group $K$. If  the initial object $\emptyset \in \cat{S}$ is removed, then the full subcategory obtained can be naturally associated with a simple graph $(V_\cat{S}, E_\cat{S})$.
Let the vertices $V_\cat{S}$ be given by $\{ 1, \ldots, r \}$ where $r$ is the rank of the maximal torus
and the edges $E_\cat{S}$ be the cardinality two elements of $\cat{S}$.
If, more generally, $K$ is essentially 2-spherical Kac--Moody group $K$, then  $\overline{\cat{S}} \sset \cat{S}$ (defined in Proposition \ref{rightcofinalprop}) is right cofinal with initial object $\bullet$. The full subcategory obtained by removing $\bullet$ can be naturally associated with a graph $(V_{\overline{\cat{S}}}, E_{\overline{\cat{S}}})$.  The  vertices $V_\cat{S}$ are minimal elements of $\overline{\cat{S}}-\{ \bullet \}$
while the edges $E_{\overline{\cat{S}}}$ are maximal in $\overline{\cat{S}}-\{ \bullet \}$.

Roughly speaking, our strategy for describing $\lim^2_{\overline{\cat{S}}} M$ %for some $M:\overline{\cat{S}}^{op} \rightarrow \cat{$R$-mod}$
will be to use a maximal spanning forest for $(V_{\overline{\cat{S}}}, E_{\overline{\cat{S}}})$ in order to relate to cases covered by Lemma \ref{lim2special}.
Compare, for example, the construction of a fundamental basis for the first homology group of a simple graph \cite[4.6]{graphth}.

 \begin{lem}
 For any essentially 2-spherical Kac--Moody group with ${\overline{\cat{S}}}\sset \cat{S}$,
$M\co\overline{\cat{S}}^{op} \rightarrow \cat{$R$-mod}$,   and (edge maximal) spanning forest of $(V_{\overline{\cat{S}}}, E_{\overline{\cat{S}}})$, there is
a finite set of subposets $\cat{S}_k  \sset \overline{\cat{S}} $
which induces a canonical short exact sequence
 \beq
0 \rightarrow \ker(\eta^\ast) \rightarrowtail \lim^2_{{\overline{\cat{S}}}} M  \stackrel{\eta^\ast}{\twoheadrightarrow} \prod \lim^2_{\cat{S}_k} M \rightarrow 0
%\nonumber
\label{lim2ses}
\eeq
\nod where each $\cat{S}_k \cong \cat{C}_{n_k}$ with ${n_k}\ge 3$ and $\ker(\eta^\ast)$ can be described explicitly (see (\ref{kernelofeta}) below).
\label{lim2general}
\end{lem}
\nod{\bf Proof:}  Let us abbreviate the notation $(V_{\overline{\cat{S}}}, E_{\overline{\cat{S}}})$ to $(V,E)$ in this proof.
 Define $\mathring{E} \sset E$ to be the edges not spanned by our hypothesized forest.
For each $e_k \in \mathring{E}=\{ e_1 , \ldots, e_m\}$, there is a unique shortest
path $\gamma_k$ in the forest between the vertices of $e_k$. Define $\cat{S}_k \sset \overline{\cat{S}}$ to be the full
subposet with
\beq
\mathrm{Objects}(\cat{S}):= \{ \bullet, e_k\} \cup \{ I | I \in \gamma_k \}
\nonumber
\label{Sk}
\eeq
\nod where $I \in \gamma_k$ means that $I \in \cat{S}$ is an edge or vertex in the path $\gamma_k$.  Notice $\cat{S}_k \cong \cat{C}_{n_k}$ for $n_k$ one more than the length of $\gamma_k$ and $n_k \ge 3$ since
$(V, E)$ is simple.
Now $\eta^\ast$ is induced by the apparent  map of posets
 \beq
 \coprod_{k=1}^{m} \cat{S}_k \stackrel{\eta}{\longrightarrow} \overline{\cat{S}}.
\label{coprodinclude}
\nonumber
\eeq
\nod It remains to show that $\eta^\ast$ is surjective and compute its kernel.

\begin{figure}[!h]
\centering
\begin{tikzpicture}[scale=.8]
\def\initalrotate{0}
\def\dmtgraph#1#2{
  \draw #1 #2 +(30+\initalrotate:2)    node[anchor=\initalrotate+210]{$3$}
   --node[anchor=270+\initalrotate]{$4$}
    +(150+\initalrotate:2)
 -- node[anchor=210+\initalrotate]{$6$}
 +(270+\initalrotate:2)  node[anchor=\initalrotate+450]{$7$}
 --node[anchor=330+\initalrotate]{$10$}
  +(390+\initalrotate:2)
  --node[anchor=210+\initalrotate] {$2$}
  +(330+\initalrotate:4)node[anchor=150+\initalrotate] {$1$}
  --node[anchor=90+\initalrotate]{$10$}    +(270+\initalrotate:2)
    --node[anchor=90+\initalrotate]{$8$}     +(210+\initalrotate:4) node[anchor=\initalrotate+30]{$9$}
     --node[anchor=\initalrotate+330]{$10$}    +(150+\initalrotate:2) node[anchor=\initalrotate+330]{$5$}
    ;

    \draw[fill=black] #1 #2 +(30+\initalrotate:2) circle (.6mm)
    +(150+\initalrotate:2) circle (.6mm)
 +(270+\initalrotate:2) circle (.6mm)
+(330+\initalrotate:4) circle (.6mm)
  +(210+\initalrotate:4) circle (.6mm)
  +(0+\initalrotate:4) node[anchor=270+\initalrotate]{$0$} circle (.6mm)
    ;
  \draw [ultra thick] #1 #2 +(330+\initalrotate:4)
   -- +(30+\initalrotate:2)
    -- +(150+\initalrotate:2)
 -- +(270+\initalrotate:2)
 --  +(210+\initalrotate:4)
    ;
}

\dmtgraph{(0,0)}{(0,0)}

\end{tikzpicture}

\caption{%A discrete Morse function  on a graph.
  Labels indicate a specific discrete Morse function \cite[2.1]{DMT}
  associated with the collapse of the bold subtree.  Each vertex with Morse value less than 10 can be paired with its unique critical edge (of lesser Morse value) to iteratively collapse the bold subtree.  This corresponds to the fact that only 0, 1 and 10 are critical values (cf. \cite[2.6]{DMT}).
}
\label{pic:oDMT}
\end{figure}
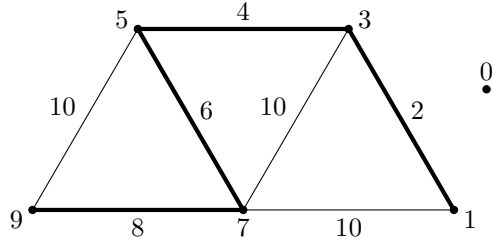

Any contraction of our hypothesized forest induces a homotopy from the metric graph realization of $(V,E)$ to a disjoint union of bouquets of circles.
In particular, we can choose a simplical collapse of the simplical realization $(V,E)$ which can be described by
(ordinary) discrete Morse theory \cite{DMT}.  See, e.g., Figure \ref{pic:oDMT}.  %For example, we could start at any end of a tree in the forest and collapse it to the adjacent $v$. Next, we can pair remaining each

%%%%%%%%%%%%%%%%%%%%%%%%%%%%%%%%%%%%%%%%%%%%%%%%%%%%%%%%%
%%%%%%%%%%%%%%%%%%%%%%%%%%%%%%%%%%%%%%%%%%%%%%%%%%%%%%%%%
%%%%%%%%%%%%%%%%%%%%%%%%%%%%%%%%%%%%%%%%%%%%%%%%%%%%%%%%%
%%%%%%%%%%%%%%%%%%%%%%%%%%%%%%%%%%%%%%%%%%%%%%%%%%%%%%%%%
%%%%%%%%%%%%%%%%%%%%%%%%%%%%%%%%%%%%%%%%%%%%%%%%%%%%%%%%%
%%%%%%%%%%%%%%%%%%%%%%%%%%%%%%%%%%%%%%%%%%%%%%%%%%%%%%%%%

Also associated to this choice is a chain homotopy equivalence from the standard normalized chain complex
\beq
{
\xymatrix{
  C^0 \ar@{=}[d] \ar[r]  & C^1 \ar[d] \ar[r]                                             & C^2  \ar[d]    \\
  C^0  \ar[r]        & \prod_{c \in \pi_0(V, E)} M_{t(c)}   \times
                                  \prod_{ E} M_v \times M_\bullet \times M_{\overline{v}} \times
                                       \prod_{\mathring{E}} M_\bullet \ar[r]                        & \prod_{\mathring{E}} M_\bullet \times M_\bullet
\label{DMT1}
} \nonumber
}
\eeq
\nod where  for an edge $e \in E$ the pair $\{v, \overline{v}\}$ is its vertices.
An example of this chain homotopy is depicted as the first stage of  Figure \ref{pic:aDMT}.  In general, this chain homotopy is induced by an acyclic pairing of $M_\bullet$ factors from $C^1$ and $C^2$ with two pairs of $M_\bullet$ factors for each edge in the spanning forest.  Our pairing is analogous to a (specific, discrete Morse theoretic)  collapse of the spanning forest.  For each pairing of $v$ and $e=\{v, \overline{v}\}$  in the (ordinary) discrete Morse function on $(V,E)$, we make two pairings, namely:
we pair the $\bullet \sset v$ coordinate with the  $\bullet \sset v \sset e$ coordinate and we pair the $\bullet \sset e$ coordinate with the  $\bullet \sset \overline{v} \sset e$ coordinate.  It is immediate from the definitions that the acyclicity of the original Morse pairing implies that this (algebraic) pairing is acyclic.

As in  Figure \ref{pic:aDMT}, we will also perform a second chain homotopy equivalence
\beq
{\xymatrix{
 C^0 \ar@{=}[d] \ar[r]  &  \prod_{c \in \pi_0(V, E)} M_{t(c)}     \times
                                     \prod_{ E} M_v   \times M_{\overline{v}}  \times
                                     \prod_{\mathring{E}} M_\bullet  \ar[r] \ar[d]                & \prod_{\mathring{E}} M_\bullet \times M_\bullet  \ar[d] \ar[r]&0   \\
  C^0  \ar[r]           &  \prod_{c \in \pi_0(V, E)} M_{t(c)}     \times
                                      \prod_{ E} M_v \times M_{\overline{v}} \ar[r]^>>>>>>>>>>>{\partial}    & \prod_{\mathring{E}} M_\bullet  \ar[r]&0
\label{DMT2}
} \nonumber
}
\eeq
\nod That is, each of the remaining $\bullet \sset \mathring{e}$ coordinates in $C^1$ for all $\mathring{e} \in \mathring{E}$ is paired with one of the two remaining $\bullet \sset  v \sset \mathring{e}$ coordinates
in $C^2$.

It remains to describe $\partial$ and its cokernel in relation to  $\cat{S}_k$.
Set $E_k$ to be the maximal elements of $\cat{S}_k$ and note that the restriction of $\partial$ to  $\prod_{c \in \pi_0(V, E)} M_{t(c)}$  is zero.
Hence, we may consider the following diagram
\beq
{\xymatrix{
 \prod_{c \in \pi_0(V, E)} M_{t(c)}     \times
\prod_{E} M_v \times M_{\overline{v}} \ar[dr]^{\partial} \ar[r]  &  \prod_{E} M_v \times M_{\overline{v}} \ar[d]^{\partial_E} \ar[r] &  \prod_{E_k} M_v \times M_{\overline{v}} \ar[d]^{\partial_{k}}       \\
                                                                      & \prod_{\mathring{E}} M_\bullet      \ar[r]             &       \prod_{\{e_k\}}      M_\bullet
}
\label{cokernel1}
 \nonumber
 }
\eeq
\nod where the unlabeled homomorphisms are the canonical projections, $\partial_E$ is defined to make the triangle commute,
and $\partial_{k}$ is associated to $\cat{S}_k$ with $\gamma_k$ chosen as a spanning tree.  The above square also
commutes because  the restriction of $\partial_E$ to the $E-E_k$ coordinates is zero in the $e_k$ coordinate. %%%% ADD FURTHER EXPLANATION
Roughly speaking, this corresponds to the fact that edges in $E-E_k$ are not incident to the 2-cell corresponding to $\cat{S}_k$.

In this way, we have identified
 $\coker(\partial_E)$ with  $\lim^2_{{\overline{\cat{S}}}} M$  and $\coker(\partial_k)$ with $\lim^2_{ \cat{S}_k} M$.
Moreover, $\{\partial_k\}_{1 \le k \le m} $ determine $\partial_E$ and with these identifications, $\eta^\ast$ corresponds to
%a quotient of the identity on $\prod_{\mathring{E}} M_\bullet$
\beq
(\prod_{\mathring{E}} M_\bullet) / \im(\partial_E)   \twoheadrightarrow (\prod_{\mathring{E}} M_\bullet) / (\prod_{\mathring{E}} \im(\partial_k))
\nonumber
\label{cokernel2}
\eeq
\nod so that
\beq
\ker(\eta^\ast) \cong \left( \prod_{\mathring{E}} \im(\partial_k) \right) / \im(\partial_E).
%\nonumber
\label{kernelofeta}
\eeq
\nod
In particular, (\ref{lim2ses}) follows and the proof is complete.
\qed

 \begin{remk} As previously noted, Lemma \ref{lim2general} depends only of the shape of $\overline{\cat{S}} $.
 For instance, Figure \ref{pic:aDMT} applies to any $D\co\cat{P}^{op} \rightarrow \cat{$R$-mod}$.
 If
\beq
\im(\partial_E)=\prod_{\{ e_1, e_2, e_3 \}} \im(\partial_k)= ( \langle D_1, D_2 , D_3  \rangle ,\langle D_2, D_3 , D_4  \rangle , \langle D_2, D_4 , D_5  \rangle ).
\nonumber
\label{examplekernelofeta}
\eeq
\nod where $D_i=D_{\{ i\}}$ as above, then $\lim^2_{\cat{P}} D  \cong \lim^2_{\cat{S}_1} D \times \lim^2_{\cat{S}_2} D \times \lim^2_{\cat{S}_3} D $.
%While we could have stated a lemma which simply gives our explicit identification of $\lim^2_{{\overline{\cat{S}}}} M$, we believe the statement of  Lemma \ref{lim2general} is
%more user friendly.

Notice also that the poset in the proof of Lemma \ref{lim2general} does not have to be finite and this
proof applies to group cohomology calculations for infinitely generated Coxeter groups which are 2-spherical (cf. Theorem \ref{nituadamsthm0}).
\qed
 \end{remk}

%%%%%%%%%%%%%%%%%%%%%%%%%%%%%%%%%%%%%%%%%%%%%%%%%%%%%%%%%
%%%%%%%%%%%%%%%%%%%%%%%%%%%%%%%%%%%%%%%%%%%%%%%%%%%%%%%%%
%%%%%%%%%%%%%%%%%%%%%%%%%%%%%%%%%%%%%%%%%%%%%%%%%%%%%%%%%
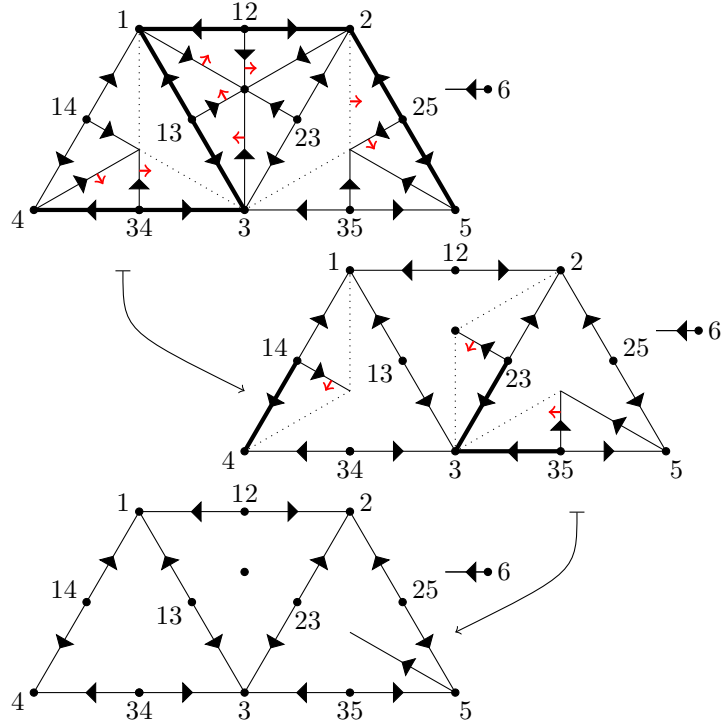
\begin{figure}[!h]
\centering
\begin{tikzpicture}[scale=.8]

\def\initalrotate{0}
\def\initalshift{(0:0)}
\def\dmtposet#1#2{
  \draw #1 #2 +(30+\initalrotate:2)    node[anchor=\initalrotate+210]{$2$}
   -- node[rotate=0+\initalrotate] {\midarrow}  +(90+\initalrotate:1) node[anchor=\initalrotate+270]{$12$}
     --node[rotate=180+\initalrotate] {\midarrow} +(150+\initalrotate:2) node[anchor=\initalrotate+330]{$1$}
  -- node[rotate=120+\initalrotate] {\midarrow} +(210+\initalrotate:1) node[anchor=\initalrotate+390]{$13$}
 -- node[rotate=300+\initalrotate] {\midarrow} +(270+\initalrotate:2)  node[anchor=\initalrotate+450]{$3$}
 -- node[rotate=240+\initalrotate]{\midarrow} +(330+\initalrotate:1) node[anchor=\initalrotate+120]{$23$}
 -- node[rotate=60+\initalrotate] {\midarrow} +(390+\initalrotate:2)

 -- node[rotate=120+\initalrotate] {\midarrow}  +(349.11+\initalrotate:2.64575) node[anchor=\initalrotate+210]{$25$}
   --node[rotate=300+\initalrotate] {\midarrow} +(330+\initalrotate:4) node[anchor=\initalrotate+120]{$5$}
 -- node[rotate=0+\initalrotate] {\midarrow} +(310.89+\initalrotate:2.64575) node[anchor=\initalrotate+90]{$35$}
 -- node[rotate=180+\initalrotate] {\midarrow}  +(270+\initalrotate:2)
  -- node[rotate=0+\initalrotate] {\midarrow}  +(229.11+\initalrotate:2.64575)node[anchor=\initalrotate+90]{$34$}
   --node[rotate=180+\initalrotate] {\midarrow} +(210+\initalrotate:4) node[anchor=\initalrotate+30]{$4$}
  -- node[rotate=240+\initalrotate] {\midarrow} +(190.89+\initalrotate:2.64575) node[anchor=\initalrotate+330]{$14$}
 -- node[rotate=60+\initalrotate] {\midarrow} +(150+\initalrotate:2)  ;

 \draw[fill=black] #1 #2 +(0+\initalrotate:4) node[anchor=180+\initalrotate]{$6$} circle (.6mm);
 \draw[fill=black] #1 #2 +(0+\initalrotate:4)  -- node[rotate=180+\initalrotate] {\midarrow} +(0+\initalrotate:3.3);
}

\def\twosimplexpic#1#2{
  \draw #1 #2 +(90+\initalrotate:2)
 -- node[rotate=120+\initalrotate] {\midarrow}  +(30+\initalrotate:1) node[anchor=\initalrotate+210]{$24$}
   --node[rotate=300+\initalrotate] {\midarrow} +(330+\initalrotate:2) node[anchor=\initalrotate+120]{$4$}
  -- node[rotate=0+\initalrotate] {\midarrow} +(270+\initalrotate:1) node[anchor=\initalrotate+90]{$13$}
 -- node[rotate=180+\initalrotate] {\midarrow} +(210+\initalrotate:2)  ;
}

\def\simplexarrowa#1#2{
  \draw #1 #2
 -- node[rotate=210+\initalrotate] {\midarrow}  +(30+\initalrotate:2)   ;
}
\def\simplexarrowb#1#2{
  \draw #1 #2
 -- node[rotate=270+\initalrotate] {\midarrow}  +(90+\initalrotate:1) ;
}

\def\simplexdottedlinea#1#2{
\draw [dotted]  #1 #2--  +(30+\initalrotate:2)   ;
}

\def\simplexdottedlineb#1#2{
\draw [dotted] #1 #2--   +(90+\initalrotate:1) ;
}

  \def\6dotunlabeledsimplex#1#2{
\draw[fill=black] #1 #2 + (30:2) circle (.6mm)
  +(90:1) circle (.6mm)
   +(150:2) circle (.6mm)
   +(210:1) circle (.6mm)
 +(270:2)  circle (.6mm)
  +(330:1) circle (.6mm) -- cycle;
  }

    \def\1dotunlabeledsimplex#1#2{
\draw[fill=black] #1 #2 +\initalshift circle (.6mm);
  }

    \def\collapse#1#2{
    \draw [ultra thick] #1 #2 +(330+\initalrotate:4)
   -- +(30+\initalrotate:2)
    -- +(150+\initalrotate:2)
 -- +(270+\initalrotate:2)
 --  +(210+\initalrotate:4)
    ;
}

    \def\collapsetwo#1#2{
    \draw [ultra thick] #1 #2 +(210+\initalrotate:4) -- +(190.89+\initalrotate:2.64575);
     \draw [ultra thick] #1 #2 +(330+\initalrotate:1) -- +(270+\initalrotate:2)
     -- +(310.89+\initalrotate:2.64575)
     ;
}

\def\simplexarrowapairing#1#2{
\draw #1 #2 [->,red, thick] ++(30+\initalrotate:.8) -- +(120+\initalrotate:.2);
}

\def\simplexarrowbpairing#1#2{
\draw [->,red, thick]  #1 #2 ++(90+\initalrotate:0.35) -- +(180+\initalrotate:0.2) ;
}

\def\simplexarrowapairingop#1#2{
\draw #1 #2 [->,red, thick] ++(30+\initalrotate:.8) -- +(300+\initalrotate:.2);
}

\def\simplexarrowbpairingop#1#2{
\draw [->,red, thick]  #1 #2 ++(90+\initalrotate:0.35) -- +(\initalrotate:0.2) ;
}

%%%%%%%%%END DEFINITIONS START DRAWING STEP ONE%%%%%%%%%%%%%%%%%%%%%%%%%%%%%%%%%%%%%%%%%
\dmtposet{(0,0)}{(0,0)}
\collapse{(0,0)}{(0,0)}

\6dotunlabeledsimplex{(0,0)}{(0,0)}

\foreach \s in {(0:0), (349.11:2.64575), (330+\initalrotate:4), (310.89+\initalrotate:2.64575),(229.11+\initalrotate:2.64575),(229.11+\initalrotate:2.64575), (210+\initalrotate:4),(190.89+\initalrotate:2.64575)} {
\def\initalshift{\s}
\1dotunlabeledsimplex{(0,0)}{(0,0)}
}

  \foreach \a in {0, 120, 240} {
  \def\initalrotate{\a}
  \simplexarrowa{(0,0)}{(0,0)}
    \simplexarrowb{(0,0)}{(0,0)}
    }

      \foreach \a in { 120, 240} {
  \def\initalrotate{\a}
  \simplexarrowapairingop{(0,0)}{(0,0)}
    }

       \foreach \a in { 0, 120} {
  \def\initalrotate{\a}
    \simplexarrowbpairingop{(0,0)}{(0,0)}
    }

      \foreach \a in {300} {
  \def\initalrotate{\a}
  \simplexarrowa{(0,0)}{(1.7320508,-1)}
  }

       \foreach \a in {60,180} {
  \def\initalrotate{\a}
    \simplexdottedlinea{(0,0)}{(1.7320508,-1)}
  }
    \foreach \a in { 60} {
  \def\initalrotate{\a}
  \simplexarrowapairingop{(0,0)}{(1.7320508,-1)}
    }

       \foreach \a in { 300} {
  \def\initalrotate{\a}
    \simplexarrowbpairingop{(0,0)}{(1.7320508,-1)}
    }

        \foreach \a in {180} {
  \def\initalrotate{\a}
  \simplexarrowa{(0,0)}{(-1.7320508,-1)}
  }

          \foreach \a in {60, 300} {
  \def\initalrotate{\a}
   \simplexdottedlinea{(0,0)}{(-1.7320508,-1)}
  }

      \foreach \a in {300, 180} {
  \def\initalrotate{\a}
  \simplexarrowb{(0,0)}{(1.7320508,-1)}
  }

        \foreach \a in {60, 180} {
  \def\initalrotate{\a}
  \simplexarrowb{(0,0)}{(-1.7320508,-1)}
  }

     \foreach \a in { 180} {
  \def\initalrotate{\a}
  \simplexarrowapairing{(0,0)}{(-1.7320508,-1)}
    }

       \foreach \a in { 180} {
  \def\initalrotate{\a}
    \simplexarrowbpairing{(0,0)}{(-1.7320508,-1)}
    }

%%%%%%%%%END STEP ONE START DRAWING STEP TWO%%%%%%%%%%%%%%%%%%%%%%%%%%%%%%%%%%%%%%%%%
\dmtposet{(0,0)}{(3.4641016,-4)}
\collapsetwo{(0,0)}{(3.4641016,-4)}
\6dotunlabeledsimplex{(0,0)}{(3.4641016,-4)}

\foreach \s in {(0:0), (349.11:2.64575), (330+\initalrotate:4), (310.89+\initalrotate:2.64575),(229.11+\initalrotate:2.64575),(229.11+\initalrotate:2.64575), (210+\initalrotate:4),(190.89+\initalrotate:2.64575)} {
\def\initalshift{\s}
\1dotunlabeledsimplex{(0,0)}{(3.4641016,-4)}
}

%%%MIDDLE
\foreach \a in {240} {
  \def\initalrotate{\a}
    \simplexarrowb{(0,0)}{(3.4641016,-4)}
    }

    \foreach \a in {0,240} {
  \def\initalrotate{\a}
  \simplexdottedlinea{(0,0)}{(3.4641016,-4)}
    }

         \foreach \a in {240} {
  \def\initalrotate{\a}
    \simplexarrowbpairingop{(0,0)}{(3.4641016,-4)}
    }

    %%RIGHT

      \foreach \a in {180} {
  \def\initalrotate{\a}
 \simplexdottedlinea{(0,0)}{(5.19615242 ,-5)}
  \simplexarrowb{(0,0)}{(5.19615242,-5)}
  }

  \foreach \a in {300} {
  \def\initalrotate{\a}
      \simplexarrowa{(0,0)}{(5.19615242,-5)}
  }

        \foreach \a in { 180} {
  \def\initalrotate{\a}
    \simplexarrowbpairingop{(0,0)}{(5.19615242,-5)}
    }

%%LEFT

        \foreach \a in {60} {
  \def\initalrotate{\a}
  \simplexarrowb{(0,0)}{(1.7320508,-5)}
  }

          \foreach \a in {60,180} {
  \def\initalrotate{\a}
   \simplexdottedlinea{(0,0)}{(1.7320508,-5)}
  }

       \foreach \a in { 60} {
  \def\initalrotate{\a}
    \simplexarrowbpairing{(0,0)}{(1.7320508,-5)}
    }

%%%%%%%%%END STEP TWO START DRAWING STEP THREE%%%%%%%%%%%%%%%%%%%%%%%%%%%%%%%%%%%%%%%%%
\dmtposet{(0,0)}{(0,-8)}
\6dotunlabeledsimplex{(0,0)}{(0,-8)}

\foreach \s in {(0:0), (349.11:2.64575), (330+\initalrotate:4), (310.89+\initalrotate:2.64575),(229.11+\initalrotate:2.64575),(229.11+\initalrotate:2.64575), (210+\initalrotate:4),(190.89+\initalrotate:2.64575)} {
\def\initalshift{\s}
\1dotunlabeledsimplex{(0,0)}{(0,-8)}
}

    %%RIGHT
  \foreach \a in {300} {
  \def\initalrotate{\a}
      \simplexarrowa{(0,0)}{(1.7320508,-9)}
  }

  %%%%%%%%%%%%%%%%%%%%%%%%%%%%%%%ARROWs

  \def\diagonalmapstoarrow#1#2{
 \draw [|->] #1 #2 .. controls +(270:1)
.. +(2,-2)  ;
}

  \def\diagonalmapstoarrowleft#1#2{
 \draw [|->] #1 #2 .. controls +(270:1)
.. +(-2,-2)  ;
}
  \diagonalmapstoarrow{(0,0)}{(-2,-3)}
  \diagonalmapstoarrowleft{(0,0)}{(5.4641016,-7)}

\end{tikzpicture}

\caption{A picture of two successive (algebraic) discrete Morse theory pairings. The top figure is the Hasse diagram for a finite poset  $\cat{P}$
isomorphic to $\overline{\cat{S}} $ covered by Lemma \ref{lim2general}. The first pairing is analogous  to the discrete Morse theory pairing in Figure \ref{pic:oDMT}
and the second pairing is a further simplification.
}
\label{pic:aDMT}
\end{figure}

 \begin{exmp}
 Let us also recall  $K:=K(A)$ from Example \ref{exmp:KA}.
Here $H^\ast(BK,  \QQ)$ is not   concentrated in odd degrees, but
\beq
H^\ast(BT,  \QQ)^W \cong \lim^0_{\cat{S}} H^\ast(BK_I,  \QQ)
\nonumber
\eeq
\nod is known to be isomorphic to $\QQ[c_2,y_4,y_6]$ additively and Theorem \ref{lim2vanish} gives that $\lim^2_{\cat{S}} H^\ast(BK_I,  \QQ) =0$.
 Hence, an Euler characteristic calculation implies that $\lim^1_{\cat{S}} H^\ast(BK_I,  \QQ)$ has Poincar\'{e} series given by:
\beq{}
\frac{1}{(1-t^2)(1-t^4)(1-t^6)}  &-&\frac{3}{(1-t^2)^2(1-t^4)(1-t^6)} + \frac{3}{(1-t^2)^3(1-t^4)} - \frac{1}{(1-t^2)^4} %-\frac{1}{(1-t^2)(1-t^4)(1-t^6)}
\nonumber \\
&=& \frac{t^4}{(1-t^2)(1-t^4)(1-t^6)}.
\nonumber
\eeq
\nod
In particular, $H^5(BK,  \QQ) \cong \lim^1_{\cat{S}} H^4(BK_I,  \QQ) $ has rank one and additively $H^\ast(BK,  \QQ)$ agrees with
$\QQ[c_2,y_4,y_6]\otimes \Lambda(x_5)$.  In fact,
it is possible to construct a weak equivalence
\beq
BK^{{\wedge}}_\QQ \simeq K(\QQ, 2) \times K(\QQ, 4) \times K(\QQ, 5) \times K(\QQ, 6);
\label{KAprod}
\eeq
\nod see Example \ref{exmp:KAprime} below.
\qed
\label{exmp:VKA}
 \end{exmp}

\section{Invariant theory and the proof of Theorem \ref{lim2vanish}}
\label{sec:invarianttheory}

The main aim of this section is to prove Theorem \ref{lim2vanish} so that, with Theorem \ref{selfmaps}, Theorem  \ref{nullconnectedintro} is proved.
Our proof of Theorem \ref{lim2vanish} will use invariant theory arguments
to show that the left and right terms of (\ref{lim2ses}) vanish for $M=H^n( D_K,  \QQ )$ with $K$  2-spherical. % Kac--Moody group .
Our primary reference for invariant theory is \cite{polyinv}, but we caution the reader that we will need to consider infinite reflection groups. %

  \begin{defn}
An involution $r \in GL_n(\QQ)$ is a reflection if its fixed space $H_r$ is $n-1$ dimensional.
  \label{reflection}
   \end{defn}

\nod
For each reflection $r$, there is a linear form $l_r\co\QQ^n \rightarrow \QQ$ is determined, up to scalar, by the fact $l_r$ has kernel $H_r$.

Let us recall how the Weyl group $W$ of a Kac--Moody group is a reflection group.
The action $W$ on the the Lie algebra $\mathfrak{h}$ of the torus is  generated by $n$ reflections
\beq
s_i(h)=h-\alpha_i(h)\alpha_i^{\vee}
\nonumber
\eeq
\nod for all $h \in \mathfrak{h}$ with  $\alpha_i \in \mathfrak{h}^\ast$ a simple root, $\alpha_i^{\vee}\in \mathfrak{h}$ a corresponding simple coroot and $n$ the size of the generalized Cartan matrix.  In particular, the simple roots \cite[1.2, 1.3]{kumar} give canonical linear forms $l_{s_i}$ and these sets of simple roots and coroots are part construction of the underlying Kac--Moody Lie algebra \cite[1.1--3]{kumar}.
For any subset $I\sset \{1, \ldots, n \}$ the subgroup $W_I:=\langle s_i | i \in I\rangle$ are the Weyl groups for corresponding parabolic
subgroups $K_I$ of $K$ which---when $W_I$ is finite---appear as $D_K(I):=BK_I$ in (\ref{hocoBK}).

Due to the following proposition, we adopt the standing convention that each $t_i$ in $\QQ[t_1, \ldots, t_r]$ has degree $2$.
 \begin{prop}
 The subgroup $W_{I}$ of the Weyl group of a Kac--Moody group is finite if and only if its collection of reflecting hyperplanes is finite.
 For any finite $W_{I}$, there is natural isomorphism $H^\ast( BK_I,  \QQ ) \cong \QQ[t_1, \ldots, t_r]^{W_I}$
 so that $H^\ast( D_K,  \QQ )$ and  $\QQ[t_1, \ldots, t_r]^{W_I}$ are naturally isomorphic functors.
\label{KacMoodyreflection}
\end{prop}
\nod{\bf Proof:}
For each reflection $r$ in $W$ the associated $l_r$ is the root associated to $r$.  In particular, none of these real roots are scalar multiples of any other \cite[1.2, 1.3]{kumar} and thus no two $H_r$ coincide.  That is, there is a one-to-one correspondence $r \leftrightarrow H_r$.

The second statement is contained in Theorem \ref{nituadamsthm0}.  In particular, the action of $W_I$ on
$H^2( BT,  \QQ )$ is dual to the action of $W_I$ on the Lie algebra $\mathfrak{h}$ of $T$ and it is well known
that the rational cohomology of the classifying space of a compact connected Lie group is given by the Weyl group invariants of $H^\ast( BT,  \QQ )$ (see, e.g., \cite{kane}).
\qed

     \begin{defn}
For a linear group $H \le GL(V)$, the relative invariants of $H$ (with respect to the determinant) $\QQ[V]^H_{\det}$ is
the set of $f\in \QQ[V]$ such that $hf= \det(h)f$ for all $h \in H$.
\label{def:relinv}
   \end{defn}
%reviewer comment
The relative invariants of $W_I$---the (possibly infinite) Weyl group of a parabolic
subgroup $K_I$ of $K$---will play a key role in our calculations.
 \begin{prop}
%Fix a field $\FF$.
 If $R \le GL_n(\QQ)$ is generated by reflections such that the collection of reflecting hyperplanes is infinite, then
 $\QQ[t_1, \ldots, t_n]^R_{\det}=0$.
\label{reflection}
\end{prop}
\nod{\bf Proof:}
%In particular, two reflections $s$ and $t$ fix the same hyperplane $H_s=H_t$ if and only if
%$s(v)=v+l_s(v)x_s$ implies $t(v)=v+l_s(v) x_t$ for some $x_t \in \QQ^n$.
%
Recall that if $h \in \QQ[t_1, \ldots, t_n]^s_{\det}$ for some reflection $s$, then
that $l_s$ divides $h$ \cite[p.226]{polyinv}. Hence, if the number of distinct reflecting hyperplanes is at least $d$, then the degree of $h$ is a least $d$.
\qed

 \begin{prop}
 If $R \le GL_n(\QQ)$ is generated by two reflections $s \neq t$ and finite, then
\beq
\QQ[t_1, \ldots, t_n]=%J \cdot \QQ[t_1, \ldots, t_n]^W
\QQ[t_1, \ldots, t_n]^R_{\det} \oplus \mathrm{span}_{\QQ}(\QQ[t_1, \ldots, t_n]^s,\QQ[t_1, \ldots, t_n]^t)
\nonumber
\eeq
\nod
as a graded vector space.
\label{reflectionsfinite}
\end{prop}
\nod{\bf Proof:}
Here $|t_i|=2$ and $R$ is a dihedral group.
  It is easy to
see that \beq
\QQ[t_1, \ldots, t_n]^R_{\det} \cap \mathrm{span}_{\QQ}(\QQ[t_1, \ldots, t_n]^s,\QQ[t_1, \ldots, t_n]^t)=0.
\eeq
\nod (cf. \cite[pp.479--480]{rank2mv}).
Hence, the internal direct sum \beq \QQ[t_1, \ldots, t_n]^R_{\det} \oplus \mathrm{span}_{\QQ}(\QQ[t_1, \ldots, t_n]^s,\QQ[t_1, \ldots, t_n]^t) \nonumber \eeq
is $\QQ[t_1, \ldots, t_n]$ if and only if these graded vector spaces have the same dimension in each degree.
By \cite[p.227]{polyinv} the Poincar\'{e} series for $\QQ[t_1, \ldots, t_n]^R_{\det}$ is %
\beq
t^{|R|} \cdot P(\QQ[t_1, \ldots, t_n]^R , t)%
\nonumber
\eeq
\nod where $P( M^\ast , t)$ denotes the Poincar\'{e} series of a graded algebra  $M^\ast$.  For any integer $k$ we have an identity of rational functions
\beq
\frac{1}{(1-x^2)(1-x^{k})}-\frac{x^{k}}{(1-x^2)(1-x^{k})} = \frac{1}{(1-x^2)} =  \frac{2}{(1-x^2)(1-x)}- \frac{1}{(1-x)^2}.
\nonumber
\eeq
\nod Considering the Poincar\'{e} series for %$\QQ[t_1, \ldots, t_n]$,
 $\QQ[t_1, \ldots, t_n]^s$, $\QQ[t_1, \ldots, t_n]^t$ and $\QQ[t_1, \ldots, t_n]^R$ (see, e.g., \cite [1.3, 7.4]{polyinv}),
the result follows from the above equality  with $x=t^2$ and $k=\frac{|R|}{2}$.
\qed

 We now summarize the invariant theory input for proving Theorem \ref{lim2vanish}.

  \begin{lem}
Let $R \le GL_n(\QQ)$ be  generated by three distinct reflections $s$, $t$, and $r$ such that
$\langle s, t\rangle$  and $\langle t, r\rangle$  are finite with $\langle s, t\rangle \cap \langle t, r\rangle=\langle t\rangle$.
  If the collection of reflecting hyperplanes in $R$ is infinite, then %for $W_{abcd}:= W_a\cap W_b \cap W_c\cap W_d$
  \beq
  \QQ[t_1, \ldots, t_n]^{\langle s, t\rangle}_{\det} \le \mathrm{span}_\QQ(\QQ[t_1, \ldots, t_n]^{r})
\label{inv1}
  \eeq
  \nod so that
  \beq
  \QQ[t_1, \ldots, t_n]=\mathrm{span}_\QQ(\QQ[t_1, \ldots, t_n]^{s},\QQ[t_1, \ldots, t_n]^{t}, \QQ[t_1, \ldots, t_n]^{r})
   \label{inv2}
  \eeq
  \nod in each degree. More generally,  if $R \le GL_n(\QQ)$ is generated by distinct reflections $r_1, \ldots, r_l, t$ for $l \ge 2$ such that
  each $\langle r_i, t\rangle$ is finite and for each pair $i \neq j$, $1\le i,j\le l$ the
  collection of reflecting hyperplanes in $\langle t, r_i, r_j\rangle$ is infinite and  $\langle r_i, t\rangle \cap \langle t, r_j\rangle=\langle t\rangle$, then
  \beq
\bigoplus_{1 \le k \le l, k\neq i } \QQ[t_1, \ldots, t_n]^{\langle r_k, t \rangle}_{\det}  \le \QQ[t_1, \ldots, t_n]^{r_i}
\label{inv3}
  \eeq
  \nod where $\oplus$ denotes internal direct sum.
  \label{lim2reflect}
   \end{lem}
\nod {\bf Proof:} By Proposition \ref{reflectionsfinite}, $\QQ[t_1, \ldots, t_n]$  coincides with
  \beq
 \QQ[t_1, \ldots, t_n]^{\langle s, t\rangle}_{\det} \oplus \mathrm{span}_\QQ(\QQ[t_1, \ldots, t_n]^s,\QQ[t_1, \ldots, t_n]^t)
  \nonumber
  \eeq
\nod in each degree  and, in particular, induces a spliting of $\mathrm{span}_\QQ(\QQ[t_1, \ldots, t_n]^{r})$.  Considering the splitting
\beq
 \QQ[t_1, \ldots, t_n]^{\langle r, t\rangle}_{\det} \oplus \mathrm{span}_\QQ(\QQ[t_1, \ldots, t_n]^r,\QQ[t_1, \ldots, t_n]^t)
  \nonumber
  \eeq
\nod (\ref{inv1}) is equivalent to
\beq
 \QQ[t_1, \ldots, t_n]^{\langle s, t\rangle}_{\det} \cap \QQ[t_1, \ldots, t_n]^{\langle r, t\rangle}_{\det}=
 \QQ[t_1, \ldots, t_n]^{R}_{\det}
= 0
  \nonumber
  \eeq
\nod
which follows from Proposition \ref{reflection}. Thus, we have (\ref{inv1}) and (\ref{inv2}) is immediate.

We prove (\ref{inv3}) by induction on $l$ and have just discussed the $l=2$ base case.  Now assume
   that for $i \neq j$, $1\le i,j\le l$  we have
    \beq
\bigoplus_{1 \le k \le l, j\neq k\neq i } \QQ[t_1, \ldots, t_n]^{\langle r_k, t \rangle}_{\det}  \le \QQ[t_1, \ldots, t_n]^{r_i}.
\label{inv4}
  \eeq
  \nod
   By the base case, we have
   \beq
  \QQ[t_1, \ldots, t_n]^{\langle r_j, t \rangle}_{\det} + \bigoplus_{1 \le k \le l, j\neq k\neq i } \QQ[t_1, \ldots, t_n]^{\langle r_k, t \rangle}_{\det}  \le \QQ[t_1, \ldots, t_n]^{r_i}
  \nonumber
    \eeq
    \nod for any $1\le i \le l$. To see that
     \beq
   \QQ[t_1, \ldots, t_n]^{\langle r_j, t \rangle}_{\det} \cap \bigoplus_{1 \le k \le l, j\neq k\neq i } \QQ[t_1, \ldots, t_n]^{\langle r_k, t \rangle}_{\det}  = 0
  \label{direct}
    \eeq
    as graded vector spaces, we note that Proposition \ref{reflectionsfinite} implies
     \beq
\QQ[t_1, \ldots, t_n]^{\langle r_j, t \rangle}_{\det} \cap \QQ[t_1, \ldots, t_n]^{r_j} = 0
  \nonumber
    \eeq
    \nod so that (\ref{direct}) follows from (\ref{inv4}) with $i$ and $j$ interchanged.
\qed

Recall that Example \ref{2sphericalcont} describes  how the somewhat technical looking statement of Lemma \ref{lim2reflect}
applies to a concrete situation.  Notice also that in Example \ref{2sphericalcont}  only (\ref{inv2}) is needed to obtain the desired
vanishing.   We will use (\ref{inv3}) to handle cases where $|E_{\cat{S}}| > > |V_{\cat{S}}|$ for the graph $(V_{\cat{S}}, E_{\cat{S}})$.
Lemmas \ref{lim2special}, \ref{lim2general} and  \ref{lim2reflect} now combine to prove Theorem \ref{lim2vanish}.

\nod {\bf Proof of Theorem \ref{lim2vanish}:}   By Proposition \ref{KacMoodyreflection}, $\lim^2_\cat{S} H^n( D_K,  \QQ )=0$ for all $n\ge 0$ if and only if
$\lim^2_\cat{S} \QQ[t_1, \ldots, t_r]^{W_I}=0$.
 For  $K$ a  2-spherical Kac--Moody group, fix a
spanning forest for  $(V_{\cat{S}}, E_{\cat{S}})$. Construct $\cat{S}_k$ as in Lemma \ref{lim2general} to obtain (\ref{lim2ses}) for the functor $M$ defined via $I \mapsto \QQ[t_1, \ldots, t_r]^{W_I}$.

Thus, $\lim^2 H^n( D_K,  \QQ )=0$ if and only if
\beq
\im(\partial_E) = \prod_{\mathring{E}} \im(\partial_k)= \prod_{\mathring{E}} M_\emptyset.
%\nonumber
\label{tight}
\eeq
\nod where here and below we use the definitions from the  proof of Lemma \ref{lim2general} freely.

% %
 Set $I_k=\bigcup_{I \in \cat{S}_k} I$  and note that $W_{I_k}:=\langle s_i | i \in I_k \rangle$ is infinite since $|I_k|=n_k \ge 3$ and $K$ is 2-spherical.
 Indeed, Lemma \ref{lim2reflect} applies to some $W_I$ for $I \sset I_k$ with $|I|=3$ so that the right equality of (\ref{tight})  now follows easily from (\ref{inv2star}) and  (\ref{inv2}) given in Lemmas \ref{lim2special} and \ref{lim2reflect}, respectively.
Define $V' \sset V_{\cat{S}}$ to be the  set of vertices not incident with some edge in $\mathring{E}$  and note
 \beq
\prod_{\mathring{E}} M_\emptyset / \im(\partial_E) &\cong& \left( \prod_{ e_k \in \mathring{E}} M_\emptyset/\partial_E(M_{v_k} + M_{{\overline{v}}_k}) \right)  /
                                                        \left(\im(\partial_E)  / \prod_{ e_k \in \mathring{E}}\partial_E(M_{v_k} + M_{{\overline{v}}_k})\right)  \nonumber\\
                                                    &\cong& \prod_{ I \in \mathring{E}} \QQ[t_1, \ldots, t_n]^{W_I}_{\det} / \partial_E\left( \prod_{ e \supseteq \{i \} \in V' } \QQ[t_1, \ldots, t_n]^{s_i}\right)
%\nonumber
\label{tight2}
\eeq
\nod where the bottom right product is taken over pairs in $\{ (e,v) | v \sset e \} \sset  E \times V'$ and $s_i$ generates the Weyl group of $K_{v=\{i \}}$. %, i.e. there is one term in the product  
The vanishing of (\ref{tight2}) is equivalent to the left equality of (\ref{tight}).

Now let $v$  be an arbitrary  leaf in the spanning forest, i.e. $v$ is in the forest and adjacent to a unique vertex in the forest. In other words, $v$ is incident with some edge in $\mathring{E}$ or
 $v$ is adjacent to a unique vertex in $(V_{\cat{S}}, E_{\cat{S}})$.
Define $L_v:=\{ v_0, \ldots, v_{|L_v|-1} \}$ to be vertices adjacent to $v$ where $v_0$ is the unique vertex in the forest adjacent to $v$. See Figure \ref{pic:aDMT}.

By design, Lemma \ref{lim2reflect} implies
\beq
 \bigoplus_{L_v-\{v_0\}} \QQ[t_1, \ldots, t_n]^{\langle s_{v_j}, s_v \rangle}_{\det}  \le \QQ[t_1, \ldots, t_n]^{s_{v_0}}
 \label{express}
 \eeq
\nod
In particular, we can partition $\mathring{E}= \coprod \mathring{E}_v$ into parts $\mathring{E}_v$ such that  all edges in $\mathring{E}_v$ adjacent to a leaf $v$. Since $v_0 \in V'$, the expression (\ref{express}) implies
that $\partial_E$ maps $\QQ[t_1, \ldots, t_n]^{s_{v_0}}$ onto $\prod_{ I \in \mathring{E}_v} \QQ[t_1, \ldots, t_n]^{W_I}_{\det}$.
Hence $\partial_E$ is surjective so that (\ref{tight2}) vanishes as desired, the left hand equality of (\ref{tight}) holds, and the proof is complete.
\qed

%%%%%%%%%%%%%%%%%%%%%%%%%%%%%%%%%%%%%%%%%%%%%%%%%%%%%%%%%
\begin{figure}[!h]
\centering
\begin{tikzpicture}[scale=.8]

\def\simplexa#1#2{
  \draw #1 #2--   +(30+\initalrotate:2)   ;
}

  \def\4dot#1#2{
\draw[fill=black] #1 #2 + (30:2) circle (.6mm) node[anchor=180]{$v_{|L_v|-1}$}
+ (126:2) circle (.6mm) node[anchor=180]{$v_2$}
   +(150:2) circle (.6mm) node[anchor=0]{$v_1 \notin V'$}
 +(270:2)  circle (.6mm) node[anchor=180]{$v_0 \in V'$}
  -- cycle;
  }
    \def\onesmalldot#1#2{
\draw[fill=black] #1 #2 + (30+\initalrotate:1.4) circle (.2mm);
  }
      \def\onetinylabeleddot#1#2{
\draw[fill=black] #1 #2 + (30+\initalrotate:1.2) circle (.01mm) node[anchor=60]{$e \in \mathring{E}$} ;
  }

        \def\onetinylabeleddottwo#1#2{
\draw[fill=black] #1 #2 + (270:0.9) circle (.01mm) node[anchor=180]{$e' \notin \mathring{E}$} ;
  }

    \def\1dotlabeled#1#2{
\draw[fill=black] #1 #2 +(0:0) circle (.6mm) node[anchor=150]{$v$};
  }

      \foreach \a in {0,96, 120, 240} {
  \def\initalrotate{\a}
  \simplexa{(0,0)}{(0,0)}
  }

       \foreach \a in {120} {
  \def\initalrotate{\a}
  \onetinylabeleddot{(0,0)}{(0,0)}
  }
     \foreach \a in {24, 36, 48, 60, 72} {
  \def\initalrotate{\a}
 \onesmalldot{(0,0)}{(0,0)}
  }

  \onetinylabeleddottwo{(0,0)}{(0,0)}

\4dot{(0,0)}{(0,0)}
\1dotlabeled{(0,0)}{(0,0)}

 \end{tikzpicture}

\caption{A local picture of $(V_{\cat{S}}, E_{\cat{S}})$ around  some leaf $v$ in the spanning forest.
}
\label{pic:aDMT}
\end{figure}
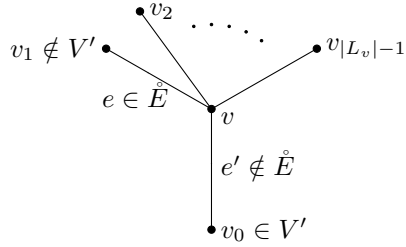

As noted before the  proof of Theorem \ref{lim2vanish}, the expression (\ref{inv3}) is particularly strong; for instance, it is straightforward to exhibit 2-spherical Kac--Moody groups so that the associated graph $(V_{\cat{S}}, E_{\cat{S}})$ is complete on $k$ vertices for any $k\in \{3, 4, 5, \ldots \}$. There is little reason to expect that $\lim^2_{\cat{S}} H^\ast( W_I ,  \QQ )$ will vanish for a general Coxeter  group  $W \le GL_n(\QQ)$ whose abstract reflections as a Coxeter group are concrete reflections in the sense of  Definition  \ref{reflection}.
 However, one may wonder what analogs of (\ref{inv3}) % and (\ref{inv3strong})
   exist for such $W$ whose abstract hyperplanes (represented, e.g., by walls in the Coxeter complex \cite{Davis}) are in bijection
   with the concrete hyperplanes $H_r$. % associated to reflections in the sense of Definition  \ref{reflection}.
   For the Weyl group of a Kac--Moody group, we pose the following question.

  \begin{quest}
For which Kac--Moody groups $K$ does $\lim^2 H^\ast( D_K,  \QQ )$ vanish? More generally, when does $\lim^{2n} H^\ast( D_K,  \QQ )=0$ for $n \ge 1$?
  \label{quest:limvanish}
   \end{quest}

\nod Notice that for any essentially 3-spherical $K$ a positive answer to the first question  implies a positive answer to the second.   Notice also that if $\lim^{2n} H^\ast( D_K,  \QQ )$ vanishes for all $n \ge 1$, then all nilpotent elements of $H^\ast( BK,  \QQ )$ appear in odd degrees so that  all products of nilpotent elements of $H^\ast( BK,  \QQ )$ vanish and the $H^\ast( BT_K,  \QQ )^W$--module structure of $H^\ast( BK,  \QQ )$
determines $H^\ast( BK,  \QQ )$ as a ring (cf. Example \ref{exmp:KAprime} below).
The construction of examples which are not essentially 3-spherical begins with $n \times n$ generalized Cartan matrices such that $n \ge 5$.

 \begin{ack}
 %Special thanks are due to Jesper Grodal.
 We are thankful to %Nitu Kitchloo,
Nat\`{a}lia Castellana,
 R\'{e}mi Molinier, Tomasz Prytula,  Albert Ruiz,  and especially Jesper Grodal for interesting discussions about this project.
Support by the Danish National Research Foundation (DNRF) through the Centre for Symmetry and Deformation (DNRF92) is gratefully acknowledged.
 \end{ack}

 \appendix
 \section{Rational product decompositions and Kac--Moody groups associated to derived Lie algebras}
 \label{append}

 As stated in the introduction, this paper studies connected unitary Kac--Moody groups as constructed in \cite[7.4]{kumar}. However, it is also common to consider (simply connected) unitary Kac--Moody groups $K'(A)$ associated to derived Kac--Moody Lie algebras (see, e.g. \cite[\S 2.5]{kacintegrate}).  Some of our references---notably \cites{nituthesis,chow,ratcohoWeyl}---follow this practice.
 This appendix will clarify  how our stated results and examples apply to this alternative setting, how  \cites{chow,ratcohoWeyl} apply to our setting, and why the rationalization of both types of groups split as a product of $K(\pi,n)$'s.   Throughout this appendix, $K$ will always denote a unitary Kac--Moody group as in \cite[7.4]{kumar} with standard torus $T$ while $K'$ will always denote a simply connected unitary Kac--Moody group as in \cite[\S 2.5]{kacintegrate} with standard torus $T'$.

 For $K$ associated to a Kac--Moody Lie algebra $\mathfrak{g}$ such that $K'$ is associated to the derived Kac--Moody Lie algebra $\mathfrak{g}':=[\mathfrak{g},\mathfrak{g}]$,
 the inclusion $K' \le K$ induces a homeomorphism of $CW$--complexes $K'/T' \cong K/T$ (cf. \cite[Appendix]{nitutkm},\cite[1.13]{nituthesis}).
 In fact, there are corresponding inclusions $K_I' \le K_I$ of compact groups for all $I \in \cat{S}$ inducing homeomorphisms $K'/K_I' \cong K/K_I$ which
 give a homeomorphism of Tits buildings
  \beq  \hcl_{I\in \cat{S}} K'/K_I' \stackrel{\approx}{\longrightarrow}  \hcl_{I\in \cat{S}} K/K_I
 \nonumber
 \eeq
 \nod and the quotient $\overline{T}: = T /T'$ is a torus of rank equal to corank of the generalized Cartan matrix  so that $T \cong T' \times \overline{T}$ (cf. \cite[5.1]{timthesis}).

    \begin{remk}  Since each $(BK')\pcp$ can be constructed as in (\ref{BXdefined}) (cf. \cite[4.2.4]{nituthesis} noting \cite[2.9]{PhDart})
    and  (\ref{BXtori}) holds with respect to $BT'$,  all the results of Sections \ref{sec:induct}-\ref{sec:nullkacmoody} apply to the alternative setting.   In particular, Theorem \ref{intronullconnectedBKp} holds with either or both of $K$ and $\Y$ replaced with
    some $K'$ as in \cite[\S 2.5]{kacintegrate}.
\qed
 \end{remk}

 Because $K$ is a finite-type \cite[Appendix]{nitutkm} $H$-space \cite{kane}, its fundamental group is finitely generated abelian  and its covering spaces are well-understood \cite[I: \S 3]{kane}.
Path lifting gives any cover of $K$, and in particular the universal cover $\widetilde{K}$, the structure of a topological group \cite{covergroups}. % REF?
Moreover, the map of fiber sequences
\beq
{\xymatrix{
T' \ar@{>->}[r] \ar[d] & T \ar[d]  \\
K' \ar@{>->}[r] \ar[d] & K \ar[d]  \\
K'/T' \ar[r]^{\approx}  & K/T
\label{triples}
}
}
\eeq
\nod yields the identification $\pi_1(\overline{T}) \cong \pi_1(K)$.

 \begin{prop}  With the preceding notation, any $K$ or $K'$ is rationally a countable product
$\prod_{n>0} K(V_n, n)$
with each $V_n$ a finite dimensional rational vector space with basis $J_n$.
\label{prop:prod}
 \end{prop}
  \nod{\bf Proof:} Since $K'$ is simply-connected and $H^\ast(K', \QQ)$ is freely generated as a ring,  the algebra generators induce a natural map
  $K' \rightarrow\prod_{n>1} K(V_n, n)$ which induces a weak equivalence $(K')\ratcp \simeq \prod_{n>1} K(V_n, n)$ \cite{rathomotopy}.

     By
\cite[I: \S3-2]{kane}, there is an equivalence $\widetilde{K} \times \overline{T} \simeq K$ so that $K$ also has this form with $|J_1|=\rank(\overline{T})$.%
  \qed

\nod Of course, Proposition \ref{prop:prod} applies to any finite type $H$-space.

 \begin{prop}  With the preceding notation, if $K'$ or $K$ is indefinite, indecomposable, then $J_{\rm odd}:=\bigcup_{n=2i+1>0} J_n$ is described in Table \ref{table:Jodd}.
\label{prop:J}
 \end{prop}
   \nod{\bf Proof:} In the alternative setting, these results for $K'$ appear in \cites{chow,ratcohoWeyl}.  In particular, an element in $J_3$ corresponds to a $W$--invariant bilinear form.

    For $\widetilde{K}$ the universal cover, we have a diagram of \emph{homomorphisms of topological groups} \cite{covergroups}
  \beq
{\xymatrix{
      & \widetilde{K}  \ar@{->>}[d]^p  \\
K' \ar@{>->}[r] \ar@{-->}[ur]^{\exists!} & K
} \nonumber
}
\eeq
\nod which induces ${K'} \simeq_{\QQ} \widetilde{K}$ and lets us complete the table.  More explicitly, (\ref{triples})
induces a map of short exact sequences of Hopf algebras \cite[Theorem 2.2]{nitutkm}
  \beq
\xymatrix{
 H^{\ast}(K/T, \QQ )\otimes_S \QQ \ar[r] \ar@{=}[d]  &  H^{\ast}(K, \QQ ) \ar[r] \ar[d] & \ar@{->>}[d] \Lambda(y_{1}, \ldots, y_{\rank(\overline{T})} , x_{1}, \ldots, x_{k}) \\
 H^{\ast}(K'/T', \QQ )\otimes_S \QQ \ar[r] &  H^{\ast}(K', \QQ ) \ar[r] & \Lambda(x_{1}, \ldots, x_{k})
}
\nonumber
\label{Hopfexactmap}
\eeq
\nod by construction such that the degree of each $y_{i}$ is one and $\widetilde{K} \times \overline{T} \simeq K$ implies ${K'} \simeq_{\QQ} \widetilde{K}$. In the terminology of \cite{nitutkm},
these $y_{i}$ correspond to a basis for the degree two part of the generalized invariants $\mathcal{I}$ which agrees with the degree two invariants $H^2( BT,  \QQ )^{W}$.
  \qed

   \begin{remk}  Propositions  \ref{prop:prod} and \ref{prop:J} imply that all the results of Sections \ref{sec:ratnullkacmoody}-\ref{sec:vanish} apply to the alternative setting---including Theorem \ref{nullconnectedLieintro}---with either or both of $K$ and $\Y$ replaced with some $K'$ as in \cite[\S 2.5]{kacintegrate}.
   \label{remk:56prime}
\qed
 \end{remk}

  \begin{table}
\caption{$J_{\rm odd}$ for indefinite, indecomposable $K'$ and $K$ }
\begin{center}
\begin{tabular}{ |l|l|l| }
\hline
& $K'$ & $K$  \\
\hline
symmetric & $|J_{\rm odd}|= |J_{3}|= 1$ & $|J_{3}|= 1$, $|J_{1}|= \rank(\overline{T})=|J_{\rm odd}|-1$ \\
\hline
not symmetric& $J_{\rm odd}= \emptyset$  & $|J_{\rm odd}|=|J_{1}|= \rank(\overline{T})$\\
\hline
\end{tabular}
\label{table:Jodd}
\end{center}
\end{table}

\begin{thm}
 If  $K'$ is 2-spherical, then
$\lim^2 H^\ast( D_{K'},  \QQ )=0$ where $D_{K'}$ is defined via $I \mapsto K_I'$ for $I \in \cat{S}$.
\label{lim2vanishprime}
\end{thm}
\nod{\bf Proof:} There is a canonical surjection $\lim^2 H^\ast( D_{K},  \QQ ) \twoheadrightarrow\lim^2 H^\ast( D_{K'},  \QQ )$. That is,
each $H^\ast( BT,  \QQ )^{W_I}$ is naturally a free $\QQ[J_{1}]$--module %$\mathcal{I}_2$ module
generated by $H^\ast( BT',  \QQ )^{W_I}$ for all $I\in \cat{S}$ so that $\QQ[J_{1}]$ is the kernel of this surjection.
  \qed

 \begin{remk}  In this way, Theorem \ref{lim2vanishprime}  and Theorem  \ref{selfmaps} (noting Remark \ref{remk:56prime}) imply that
  Theorem \ref{nullconnectedintro} holds with $K$ replaced by some $K'$ as in \cite[\S 2.5]{kacintegrate}.
\qed
 \end{remk}

We finish this appendix by describing how Example \ref{exmp:KA} changes when moving to the simply connected setting of \cite{kacintegrate}.  Notice that Examples \ref{1spherical} and \ref{2spherical} are not effected.
We also determine the rational type of $BK(A)$.
 \begin{exmp}
Now let $K':=K'(A)$ be the simply-connected  affine Kac--Moody group
associated to the generalized Cartan matrix
  \beq
A=\left[\begin{array}{ccc} 2 & -1 & -1 \\ -1 & 2 & -1 \\ -1 & -1 &  2 \end{array}\right] .
\nonumber
\label{gcm3}
\eeq
\nod Then
$
% \beq
 H^{\ast}( K', \QQ ) = \Lambda(y_3, y_5) \otimes \QQ[x_4]
%\nonumber
$.
Here the degree one generator in $H^{\ast}( K(A), \QQ )$ corresponds an element of $H^2( BT,  \QQ )^{W}$ and
this generator does not appear in $H^{\ast}( K'(A), \QQ )$.
With analogous definitions and Euler characteristic arguments,
\beq
\ker\left( V^{T'}_{K'} \rightarrow V^{K'}_{K'} \right) \cong H^{5}( BK', \hat{\QQ}/\QQ)\cong \hat{\QQ}/\QQ.
\eeq
\nod as in Example \ref{exmp:VKA}.  Hence, each non-zero $x \in \hat{\QQ}/\QQ$ is represented by a non-trivial $f_x\co BK' {\rightarrow} BK'$
with ${f_x}B\iota$ is null and  nontrivial
$f_x\co BK' {\rightarrow} BG$ with $G$ is compact connected Lie can be constructed as in Example \ref{exmp:KA}.

Considering the bigraded algebra structure on the
 cohomological  spectral sequence \cite[XII. 5.8]{BK}
with
\beq
F_{\ff}^{i, k} = F_2^{i, k} = \lim^i_{I \in \cat{S}} H^{k} (BK_I',  \QQ)
\label{cohomologyinexample}
\nonumber
\eeq
\nod the $H^\ast( BT,  \QQ )^W$--module structure of $H^\ast( BK'(A),  \QQ )$
determines $H^\ast( BK',  \QQ )$ as a ring.  Hence, $H^\ast( BK',  \QQ ) \cong \QQ[x_4, x_6] \otimes \Lambda(y_5)$ as a ring
so that---as in the proof of Proposition \ref{prop:prod}---we have a weak equivalence
\beq
(BK')^{{\wedge}}_\QQ \simeq K(\QQ, 4) \times K(\QQ, 5) \times K(\QQ, 6).
\eeq
\nod Similarly,  $H^\ast( BK(A),  \QQ ) \cong \QQ[x_2, x_4, x_6] \otimes \Lambda(y_5)$ as a ring and (\ref{KAprod}) follows---as in the proof of Proposition \ref{prop:prod}.
\qed
\label{exmp:KAprime}
 \end{exmp}

\end{document}